\documentclass{article}
\usepackage{amsmath}
\usepackage{amsthm}
\usepackage{amssymb}
\usepackage{bm}

\newtheorem{theorem}{Theorem}[section]
\newtheorem{proposition}[theorem]{Proposition}
\newtheorem{lemma}[theorem]{Lemma}
\newtheorem{corollary}[theorem]{Corollary}
\theoremstyle{definition}
\newtheorem{definition}[theorem]{Definition}
\theoremstyle{remark}
\newtheorem{remark}[theorem]{Remark}

\newcommand{\N}{\mathbb{N}}
\newcommand{\Z}{\mathbb{Z}}
\newcommand{\Q}{\mathbb{Q}}
\newcommand{\R}{\mathbb{R}}
\newcommand{\C}{\mathbb{C}}

\newcommand{\calD}{\mathcal{D}}
\newcommand{\calH}{\mathcal{H}}

\newcommand{\bmu}{\bm{\mu}}
\newcommand{\bnu}{\bm{\nu}}
\newcommand{\blambda}{\bm{\lambda}}
\newcommand{\balpha}{\bm{\alpha}}
\newcommand{\bbeta}{\bm{\beta}}
\newcommand{\muvec}{\mu_1,\ldots,\mu_p}
\newcommand{\nuvec}{\nu_1,\ldots,\nu_q}

\newcommand{\shp}{\,\#\,}
\newcommand{\dotshp}{\,\dot{\#}\,}
\newcommand{\barast}{\,\overline{*}\,}
\newcommand{\barcast}{\,\overline{\circledast}\,}
\newcommand{\barO}{\overline{O}}
\newcommand{\tI}{\tilde{I}}
\newcommand{\tV}{\tilde{V}}
\newcommand{\id}{\mathop{\mathrm{id}}\nolimits}
\newcommand{\End}{\mathop{\mathrm{End}}\nolimits}
\newcommand{\re}{\mathop{\mathrm{Re}}\nolimits}

\begin{document}
\title{A class of relations among multiple zeta values}
\author{Gaku Kawashima \\
        Graduate School of Mathematics, \\
        Nagoya University, Chikusa-ku, Nagoya464-8602, Japan \\
        E-mail: m02009c@math.nagoya-u.ac.jp}
\date{}
\maketitle
\begin{abstract}
    We prove a new class of relations among multiple zeta values (MZV's)
    which contains Ohno's relation. 
    We also give the formula for the maximal number of independent MZV's of fixed weight, under our new relations.
    To derive our formula for MZV's, we consider the Newton series whose
    values at non-negative integers are finite multiple harmonic sums.
\end{abstract}
\textbf{keywords} multiple zeta value, the Newton series, Ohno's relation
%
%
\section{Introduction}
Let $k_1,\ldots,k_p$ be positive integers and let $k_p \ge 2$.
Multiple zeta values (MZV's) are defined by
\begin{displaymath}
  \zeta(k_1,\ldots,k_p) = \sum_{0 < n_1 < \cdots < n_p} \frac{1}{n_1^{k_1} \cdots n_p^{k_p}}.
\end{displaymath}
The sum $k_1 + \cdots + k_p$ is called the weight of the multiple zeta value $\zeta(k_1,\ldots,k_p)$.
These numbers were considered first by Euler \cite{E}; he studied the case $p=2$.
The general case was introduced in \cite{Hof4} and \cite{Z}.
In recent years, many researchers have studied these numbers 
in connection with Galois representations, arithmetic geometry, quantum groups, invariants for knots, mathematical physics, etc.
Many $\mathbb{Q}$-linear relations are known among these values, 
e.g. Ohno's relation, the cyclic sum formula and the derivation relation \cite{HO,IKZ,Ohno1,OW}.
The famous relation due to Ohno generalizes the duality, the sum formula and Hoffman's relation simultaneously.
Ihara, Kaneko and Zagier investigated the regularized double shuffle relations in \cite{IKZ},
where they conjectured that these relations imply all
$\mathbb{Q}$-linear relations.
Recently in \cite{K} Kaneko proposed a conjectural generalization of the
derivation relation for MZV's and in \cite{T} Tanaka proved this conjecture by reducing
it to the relations studied in this paper.
The aim of this paper is to give a new class of relations among MZV's.
We shall prove that our new class of relations contains Ohno's relation.
We conjecture that it also contains the cyclic sum formula.
We have checked this up to weight 12.
\par
We shall outline how our relations among MZV's can be derived.
Let $\bm{\mu} = (\mu_1,\ldots,\mu_p)$ be a multi-index (i.e., an ordered set of positive integers) and $n$ a non-negative integer.
We consider the finite multiple harmonic sums
\begin{displaymath}
  s_{\bm{\mu}}(n) = \sum_{0 \le n_1 \le \cdots \le n_p = n} \frac{1}{(n_1+1)^{\mu_1} \cdots (n_p+1)^{\mu_p}}
\end{displaymath}
and
\begin{align*}
  S_{\bm{\mu}}(n) &= \sum_{0 \le k < n} s_{\bm{\mu}}(k) \\
             &= \sum_{0 \le n_1 \le \cdots \le n_p < n} \frac{1}{(n_1+1)^{\mu_1} \cdots (n_p+1)^{\mu_p}}.
\end{align*}
In the following, we denote the set of non-negative integers by $\mathbb{N} = \{0,1,2,\ldots \}$.
We define the operator $\nabla$ on the space $\mathbb{C}^{\mathbb{N}}$ of complex-valued sequences by
\begin{displaymath}
  (\nabla a)(n) = \sum_{k=0}^{n}(-1)^k \binom{n}{k} a(k)
\end{displaymath}
for $a \in \mathbb{C}^{\mathbb{N}}$ and $n \in \N$.
In Section \ref{multiple-sum}, we shall prove our key formula for the finite multiple harmonic sums $s_{\bm{\mu}}(n)$ (Theorem \ref{th3-10}).
As a corollary of Theorem \ref{th3-10}, we shall prove that
\begin{displaymath}
 \nabla s_{\bm{\mu}} = s_{\bmu^{*}},
\end{displaymath}
where $\bmu^{*}$ is some multi-index determined by $\bm{\mu}$.
We note that this assertion was proved independently by Hoffman in \cite[Theorem~4.2]{Hof1}.
\par
In Section \ref{multiple-zeta}, in order to derive our relations among MZV's, 
we consider the Newton series which interpolates
the finite multiple harmonic sums $S_{\bmu}(n)$. This Newton series
$F_{\bmu}(z)$ is defined by
\begin{displaymath}
  F_{\bm{\mu}}(z) = \sum_{n=0}^{\infty} (-1)^n (\nabla S_{\bm{\mu}})(n) \binom{z}{n},\quad 
                                                               \binom{z}{n} = \frac{z(z-1) \cdots (z-n+1)}{n!},
\end{displaymath}
where $z$ is a complex number with $\re z > -1$.
In fact, for any $n \in \mathbb{N}$, we have $F_{\bm{\mu}}(n) = S_{\bm{\mu}}(n)$.
It is well known that $S_{\bm{\mu}}(n) S_{\bm{\nu}}(n)$
is a $\mathbb{Z}$-linear combination of $S_{\bm{\lambda}}(n)$'s:
\begin{displaymath}
  S_{\bm{\mu}}(n) S_{\bm{\nu}}(n) = \sum_i c_i S_{\bm{\lambda}^i}(n),\quad c_i \in \mathbb{Z}.
\end{displaymath}
The multi-indices $\bm{\lambda}^i$ and the coefficients $c_i$ depend only on $\bm{\mu}$ and $\bm{\nu}$, but not on $n$.
From this fact, we derive the functional equations
\begin{equation} \label{eq1-10} 
  F_{\bm{\mu}}(z) F_{\bm{\nu}}(z) = \sum_i c_i F_{\bm{\lambda}^i}(z),
\end{equation}
which are valid for $\re z > -1$.
Since the derivatives $F_{\bm{\mu}}^{(m)}(0)$ $(m \ge 1)$ can be
expressed in terms of MZV's, quadratic relations among MZV's are derived from (\ref{eq1-10}) (Corollary \ref{cor5-10}).
\par
We mainly consider the first derivatives only.
Since $F_{\bm{\mu}}(0)$ vanishes, equation (\ref{eq1-10}) implies a linear relation
\begin{displaymath}
  0  = \sum_i c_i F_{\bm{\lambda}^i}'(0).
\end{displaymath}
In this way, we get a set of linear relations among MZV's which can be
written explicitly (Corollary \ref{cor5-20}). We shall explain this formula.
Let $V_k$ be the $\mathbb{Q}$-vector space whose basis is the set of all multi-indices of weight $k$ and let $V = \bigoplus_{k \ge 1} V_k$.
We define $\mathbb{Q}$-linear mappings $\sigma$, $u$ from $V$ to $V$ by
\begin{align*}
  \sigma(\mu_1,\ldots,\mu_p) &= (-1)^p (\mu_1,\ldots,\mu_p),\\
  u(\mu_1,\ldots,\mu_p) &= \sum_{(\nu_1,\ldots,\nu_q) \ge (\mu_1,\ldots,\mu_p)} (\nu_1,\ldots,\nu_q),
\end{align*}
where $(\nu_1,\ldots,\nu_q)$ runs over all refinements of $(\mu_1,\ldots,\mu_p)$
(see \S \ref{multi-indices} for the precise definition).
We define the $\mathbb{Q}$-linear mapping $\zeta^{+}$ from $V$ to $\mathbb{R}$ by
\begin{displaymath}
  \zeta^{+}(\mu_1,\ldots,\mu_p) = \zeta(\mu_1,\ldots,\mu_{p-1},\mu_p + 1).
\end{displaymath}
Then Corollary \ref{cor5-20} asserts that 
\begin{displaymath}
  \zeta^{+}(u \sigma(\bm{\mu} * \bm{\nu})) = 0
\end{displaymath}
for any multi-indices $\bm{\mu}$ and $\bm{\nu}$.
(The multiplication $*$ on $V$ is the harmonic product.)
For example, for $\bm{\mu} = \bm{\nu} = (1)$, we have
\begin{align*}
  \bm{\mu} * \bm{\nu} &= (2) + 2(1,1),\\
  \sigma(\bm{\mu} * \bm{\nu}) &= -(2) + 2(1,1)\\
  \intertext{and}
  u \sigma(\bm{\mu} * \bm{\nu}) &= -\{(2) + (1,1)\} + 2(1,1) = -(2) + (1,1).
\end{align*}
Therefore we obtain
\begin{displaymath}
  \zeta^{+}(2) = \zeta^{+}(1,1) \quad \bigl(\text{i.e., } \zeta(3) = \zeta(1,2) \bigr),
\end{displaymath}
which is due to Euler.
\par
Under the harmonic product $*$, the $\Q$-vector space $\Q \oplus V$
becomes a commutative $\Q$-algebra. This is isomorphic to the
commutative $\Q$-algebra of quasi-symmetric functions, which is known to
be the polynomial algebra on the set of Lyndon
words~\cite[Section~2]{Ma}, and whose enumeration is well known. From
this fact, it is easily seen that
\begin{align} \nonumber
  d(k) &:= \dim_{\mathbb{Q}}(u\sigma(V * V) \cap V_k)\\ \label{eq:intro2}
       &=  2^{k-1} - \frac{1}{k} \sum_{d|k} \mu \left( \frac{k}{d} \right) 2^d \quad (k \ge 2),
\end{align}
where $\mu$ is the M\"obius function. 
The right-hand side of (\ref{eq:intro2}) appeared in the recent research of Kaneko for the extended derivation relation,
and this equation was suggested to the author by Kaneko.
\par
We shall compare the dimensions of three subspaces of $\ker \zeta^{+} \cap V_k$.
We define the sequence $\{z_k\}_{k=1}^{\infty}$ by 
\begin{displaymath}
  z_1 = z_2 = z_3 = 1,\quad z_k = z_{k-2} + z_{k-3} \quad (k \ge 4)
\end{displaymath}
and put $d_Z(k) = 2^{k-1} - z_k$.
Then the dimension of the space $\ker \zeta^{+} \cap V_k$ is conjectured to be equal to $d_Z(k)$ by Zagier \cite{Z}.
We denote the dimension of the space $V_{Ohno} \cap V_k$ by $d_O(k)$,
where $V_{Ohno}$ is the subspace of $\ker \zeta^{+}$ which corresponds to Ohno's relation.
The following table gives us the numerical values of $d_Z(k)$, $d(k)$ and $d_O(k)$ for $2 \le k \le 11$.
\begin{center}
    \begin{tabular}{|c|r|r|r|r|r|r|r|r|r|r|} \hline
        $k$               & 2 & 3 & 4 & 5  & 6  & 7  & 8    & 9   & 10  & 11 \\ \hline
        $d_Z(k)$          & 1 & 3 & 6 & 14 & 29 & 60 & 123  & 249 & 503 & 1012\\ \hline
        $d(k)$            & 1 & 2 & 5 & 10 & 23 & 46 & 98   & 200 & 413 & 838\\ \hline
        $d_O(k)$          & 1 & 2 & 5 & 10 & 23 & 46 & 98   & 199 & 411 & 830\\ \hline
    \end{tabular}
\end{center}
\par
The product of two MZV's is expressed as a $\mathbb{Z}$-linear combination of MZV's in two ways;
the harmonic product relations and the shuffle product relations.
We conjecture that under the shuffle product relations, Corollary \ref{cor5-10} for 
$m = 1,2$ gives all $\Q$-linear relations among MZV's.
This conjecture is due to Tanaka.
He checked that the dimension of the subspace of $\ker \zeta^{+}$ which
corresponds to these linear relations 
is equal to $d_{Z}(k)$ for $2 \le k \le 11$.
%
%
\section{Multi-indices} \label{multi-indices}
A finite sequence of positive integers is called a multi-index and we
denote the set of all multi-indices by $I$.
We denote by $V$ the $\mathbb{Q}$-vector space whose basis is $I$.
For $\bm{\mu} = (\mu_1,\ldots,\mu_p) \in I$, we call 
$l(\bm{\mu}) := p$ and $|\bm{\mu}| := \sum_{i=1}^{p}\mu_i$
the length and the weight of $\bm{\mu}$, respectively. 
For a positive integer $m$, we denote the set of all multi-indices of weight $m$ by $I_m$.
We define the mapping $\mathcal{S}_m$ for a positive integer $m$ by
\begin{displaymath}
  \mathcal{S}_m \colon I_m \rightarrow 2^{\{1,2,\ldots,m-1 \}},\quad
  (\mu_1,\ldots,\mu_p) \mapsto \left\{\sum_{i=1}^k \mu_i \Bigm| 1 \le k < p \right\}.
\end{displaymath}
This is a bijection. For example, we obtain
\begin{displaymath}
  \mathcal{S}_5 (2,2,1) = \{2,4\} \quad \text{and} \quad \mathcal{S}_5 (1,1,3) = \{1,2\}
\end{displaymath}
by the following diagrams
\smallskip
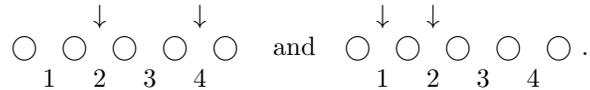
\begin{displaymath}
  {\arraycolsep=1pt
  \begin{array}{cccccccccc}
          & &        &\downarrow&        & &        &\downarrow&        &\\  
  \bigcirc& &\bigcirc&          &\bigcirc& &\bigcirc&          &\bigcirc&\\
          &1&        &2         &        &3&        &4         &        &
  \end{array}}
  \quad \text{and} \quad
  {\arraycolsep=1pt
  \begin{array}{cccccccccc}
          &\downarrow&        &\downarrow&        & &        & &        &\\
  \bigcirc&          &\bigcirc&          &\bigcirc& &\bigcirc& &\bigcirc&\,.\\
          &1         &        &2         &        &3&        &4&        &
  \end{array}}
\end{displaymath}
\smallskip
We define an order on $I$ by setting $\bm{\mu} \ge \bm{\nu}$ 
if $|\bm{\mu}| = |\bm{\nu}|$ and $\mathcal{S}_{|\bm{\mu}|}(\bm{\mu}) \supset \mathcal{S}_{|\bm{\nu}|}(\bm{\nu})$.
In other words, we write $\bm{\mu} \ge \bm{\nu}$ if $\bm{\mu}$ is a refinement of $\bm{\nu}$.
For example, the order on $I_4$ is given by the following Hasse diagram
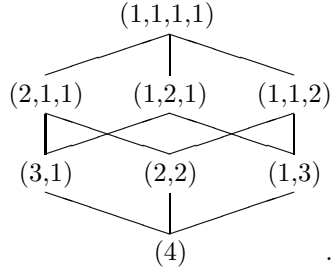
\begin{center}
  \begin{picture}(120,100)
    \put(55,0){(4)}
    \put(4,30){(3,1)}
    \put(51,30){(2,2)}
    \put(98,30){(1,3)}
    \put(0,60){(2,1,1)}
    \put(47,60){(1,2,1)}
    \put(94,60){(1,1,2)}
    \put(42,90){(1,1,1,1)}
    \put(61,10){\line(0,1){15}}
    \put(61,70){\line(0,1){15}}
    \put(14,40){\line(0,1){15}}
    \put(108,40){\line(0,1){15}}
    \put(61,10){\line(-3,1){47}}
    \put(61,10){\line(3,1){47}}
    \put(61,40){\line(-3,1){47}}
    \put(61,40){\line(3,1){47}}
    \put(14,40){\line(3,1){47}}
    \put(108,40){\line(-3,1){47}}
    \put(14,70){\line(3,1){47}}
    \put(108,70){\line(-3,1){47}}
    \put(120,0){.}
  \end{picture}
\end{center}
For $\bm{\mu} = (\mu_1,\ldots,\mu_p) \in I$ with $|\bm{\mu}| = m$, we put
\begin{align*}
  *(\bm{\mu}) &= \bm{\mu}^{*} = \mathcal{S}_{m}^{-1}(\mathcal{S}_{m}(\bm{\mu})^c),\\
  \tau(\bm{\mu}) &= (\mu_p,\ldots,\mu_1),\\
  \sigma(\bm{\mu}) &= (-1)^{l(\bm{\mu})}\bm{\mu},\\
  u(\bm{\mu}) &= \sum_{\bm{\nu} \ge \bm{\mu}}\bm{\nu}\\
\intertext{and}
  d(\bm{\mu}) &= \sum_{\bm{\nu} \le \bm{\mu}}\bm{\nu}.
\end{align*}
($\mathcal{S}_{m}(\bm{\mu})^c$ denotes the complement of $\mathcal{S}_{m}(\bm{\mu})$ in $\{1,2,\ldots,m-1\}$.)
We extend them as $\mathbb{Q}$-linear mappings from $V$ to $V$.
All of these are bijections. 
We give examples of $\bmu^{*}$.
We have
\begin{displaymath}
  (1,2,3)^{*} = (2,2,1,1),\quad (2,2,2)^{*} = (1,2,2,1)\quad \text{and}\quad (4,1,1)^{*} = (1,1,1,3)
\end{displaymath}
by the diagrams
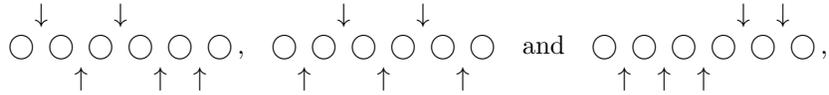
\begin{displaymath}
 {\setlength{\arraycolsep}{0pt}
 \begin{array}{cccccccccccccccccccccccccccccccccccc}
  &\downarrow& & & &\downarrow& & & & & & & & & &\downarrow& & & &\downarrow& & & & & & & & & & & &\downarrow& &\downarrow& &
   \\
  \bigcirc& &\bigcirc& &\bigcirc& &\bigcirc& &\bigcirc& &\bigcirc&\, ,\quad&\bigcirc& &\bigcirc& &\bigcirc& &\bigcirc& &\bigcirc& &\bigcirc&\quad\text{and}\quad&\bigcirc& &\bigcirc& &\bigcirc& &\bigcirc& &\bigcirc& &\bigcirc&\,, \\
  & & &\uparrow& & & &\uparrow& &\uparrow& & & &\uparrow& & & &\uparrow& & & &\uparrow& & & &\uparrow& &\uparrow& &\uparrow& & & & & & \\
 \end{array}}
\end{displaymath}
where the lower arrows are in the complementary slots to the upper arrows.
We can also calculate these correspondences by the diagrams
\smallskip
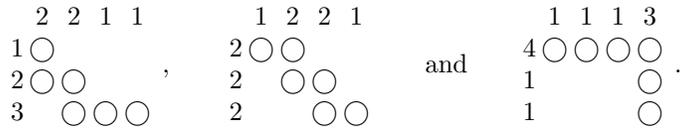
\begin{displaymath}
  {\arraycolsep=1pt
  \begin{array}{ccccc}
     &2       &2       &1       &1        \\
    1&\bigcirc&        &        &         \\
    2&\bigcirc&\bigcirc&        &         \\
    3&        &\bigcirc&\bigcirc&\bigcirc
  \end{array}}\,\, ,\quad\quad
  {\arraycolsep=1pt
  \begin{array}{ccccc}
     &1       &2       &2       &1       \\
    2&\bigcirc&\bigcirc&        &        \\
    2&        &\bigcirc&\bigcirc&        \\
    2&        &        &\bigcirc&\bigcirc
  \end{array}} \quad\quad \text{and} \quad\quad
  {\arraycolsep=1pt
  \begin{array}{ccccc}
     &1       &1       &1       &3        \\
    4&\bigcirc&\bigcirc&\bigcirc&\bigcirc \\
    1&        &        &        &\bigcirc \\
    1&        &        &        &\bigcirc
  \end{array}}\,\, .
\end{displaymath}
\smallskip
It is easily seen that
\begin{equation}
 l(\bmu) + l(\bmu^{*}) = |\bmu| + 1 \label{eq2-5} 
\end{equation}
for any $\bmu \in I$.
Let
\begin{displaymath}
 \bmu = (\mu_1,\ldots,\mu_r,\underbrace{1,\ldots,1}_l) \quad (\mu_r \ge 2,\,
 r \ge 0,\, l \ge 0)
\end{displaymath}
and $\bmu^{*} = (\mu^{*}_1,\ldots,\mu^{*}_q)$. Then we note that
\begin{equation} \label{eq2-10} 
 \mu^{*}_q = l+1.
\end{equation}
\begin{proposition}  
 \label{eq2-15} 
 We have $*\tau= \tau *$, $\sigma \tau = \tau \sigma$,
 $u \tau = \tau u$ and $d \tau = \tau d$.
\end{proposition}
\begin{proof}
 These are easily seen.
\end{proof}
\begin{lemma} \label{lem2-10}
  Let $T_1 \subset T_2$ be finite sets. Then we have
  \begin{displaymath}
    \sum_{T_1 \subset S \subset T_2} (-1)^{\# S} = \begin{cases}
                                                     (-1)^{\# T_2} & \text{if $T_1 = T_2$} \\
                                                     0             & \text{if $T_1 \subsetneq T_2$}.
                                                   \end{cases}
  \end{displaymath}
\end{lemma}
\begin{proof}
  The left-hand side equals 
  \begin{align*}
    (-1)^{\# T_1} \sum_{S \subset (T_2 \setminus T_1)} (-1)^{\# S} 
    &= (-1)^{\# T_1} \sum_{k=0}^{\# (T_2 \setminus T_1)} (-1)^k \binom{\# (T_2 \setminus T_1)}{k} \\
    &= (-1)^{\# T_1} (1-1)^{\# (T_2 \setminus T_1)}.
  \end{align*}
  Thus we proved the lemma.
\end{proof} 
\begin{proposition} \label{prop2-20} 
  We have\\
  \rm{(i)} $*d*=u$,\\
  \rm{(ii)} $d \sigma d \sigma =\id_V$,\\
  \rm{(iii)} $u \sigma u \sigma =\id_V$,\\
  \rm{(iv)} $d*d^{-1}=-u \sigma$,\\
  \rm{(v)} $u^{-1}*u=- \sigma d$.
\end{proposition}
\begin{proof}
  Let $\bm{\mu}$ be any element of $I$ and let $|\bm{\mu}| = m$.\\
  (i) We have
  \begin{displaymath}
    *d*(\bm{\mu}) = \sum_{\bm{\nu} \le \bm{\mu}^{*}} \bm{\nu}^{*} 
                  = \sum_{\bm{\nu}^{*} \ge \bm{\mu}} \bm{\nu}^{*}
                  = u(\bm{\mu}).
  \end{displaymath}
  (ii)We have 
  \begin{displaymath}
    d \sigma d \sigma (\bm{\mu}) 
    = (-1)^{l(\bm{\mu})} \sum_{\bm{\nu} \le \bm{\mu}} (-1)^{l(\bm{\nu})} \sum_{\bm{\lambda} \le \bm{\nu}} \bm{\lambda} 
    = (-1)^{l(\bm{\mu})} \sum_{\bm{\lambda} \le \bm{\mu}} 
       \left\{\sum_{\bm{\lambda} \le \bm{\nu} \le \bm{\mu}} (-1)^{l(\bm{\nu})} \right\} \bm{\lambda}.
  \end{displaymath}
  By Lemma \ref{lem2-10}, we have
  \begin{displaymath}
    \sum_{\bm{\lambda} \le \bm{\nu} \le \bm{\mu}} (-1)^{l(\bm{\nu})}
    = \sum_{\mathcal{S}_m(\bm{\lambda}) \subset \mathcal{S}_m(\bm{\nu}) \subset \mathcal{S}_m(\bm{\mu})} (-1)^{\# \mathcal{S}_m(\bm{\nu}) + 1}
    = \begin{cases}
        (-1)^{l(\bm{\mu})} & \text{if $\bm{\lambda} = \bm{\mu}$} \\
        0                  & \text{if $\bm{\lambda} < \bm{\mu}$}.
      \end{cases}
  \end{displaymath}
  Therefore we obtain $d \sigma d \sigma (\bm{\mu}) = \bm{\mu}$. \\
  (iii) The proof of (iii) is similar to the proof of (ii). \\
  (iv) We have
  \begin{displaymath}
    d* \sigma d(\bm{\mu})
    = \sum_{\bm{\nu} \le \bm{\mu}} (-1)^{l(\bm{\nu})} \sum_{\bm{\lambda} \le \bm{\nu}^{*}} \bm{\lambda} 
    = \sum_{\bm{\lambda} \in I_m} 
    \left\{ \sum_{\substack{\bm{\nu} \le \bm{\mu}\\ \bm{\nu} \le \bm{\lambda}^{*}}} (-1)^{l(\bm{\nu})} \right\} \bm{\lambda}.
  \end{displaymath}
  By Lemma \ref{lem2-10}, we have
  \begin{displaymath}
    \sum_{\substack{\bm{\nu} \le \bm{\mu}\\ \bm{\nu} \le \bm{\lambda}^{*}}} (-1)^{l(\bm{\nu})}
    = \sum_{\mathcal{S}_m(\bm{\nu}) \, \subset \, \mathcal{S}_m(\bm{\mu}) \cap \mathcal{S}_m(\bm{\lambda})^{c}} (-1)^{\# \mathcal{S}_m(\bm{\nu}) + 1}
    = \begin{cases}
        -1 & \text{if $\mathcal{S}_m(\bm{\mu}) \subset \mathcal{S}_m(\bm{\lambda})$} \\
        0  & \text{if $\mathcal{S}_m(\bm{\mu}) \not\subset \mathcal{S}_m(\bm{\lambda})$}.
      \end{cases}
  \end{displaymath}
  Hence we get
  \begin{displaymath}
    d* \sigma d(\bm{\mu}) = -\sum_{\bm{\lambda} \ge \bm{\mu}} \bm{\lambda} = -u(\bm{\mu}).
  \end{displaymath}
  This implies the assertion.\\
  (v) The proof of (v) is similar to the proof of (iv).
\end{proof}
We define the multiplication $*$ on $V$ known as the harmonic product.
We shall give a combinatorial description of the harmonic product.
For $\bm{\mu}$, $\bm{\nu} \in I$, we consider $2 \times l$ $(l \ge 1)$ matrices $M = (m_{ij})$ with the following properties:
\begin{itemize}
  \item Each entry of $M$ is a non-negative integer.
  \item If we exclude zeros from the first row of $M$, we get $\bm{\mu}$.
  \item If we exclude zeros from the second row of $M$, we get $\bm{\nu}$.
  \item For $1 \le j \le l$, we have $m_{1j} + m_{2j} > 0$.
\end{itemize}
(Then it must hold that $\max\{l(\bm{\mu}),l(\bm{\nu})\} \le l \le l(\bm{\mu}) + l(\bm{\nu})$.)
We denote the set of all such matrices by $\mathcal{H}_{\bm{\mu}, \bm{\nu}}$.
For example, for $\bm{\mu} = (1)$ and $\bm{\nu} = (2,3)$, we have
\begin{equation}
  \mathcal{H}_{\bm{\mu}, \bm{\nu}}
  = \left\{ \begin{pmatrix} 1 & 0 & 0 \\ 0 & 2 & 3 \end{pmatrix},\,
            \begin{pmatrix} 0 & 1 & 0 \\ 2 & 0 & 3 \end{pmatrix},\,
            \begin{pmatrix} 0 & 0 & 1 \\ 2 & 3 & 0 \end{pmatrix},\,
            \begin{pmatrix} 1 & 0 \\ 2 & 3 \end{pmatrix},\,
            \begin{pmatrix} 0 & 1 \\ 2 & 3 \end{pmatrix} 
  \right\}. \label{eq2-25}
\end{equation}
\begin{definition} 
 \label{def2-30}
  For $\bm{\mu}$, $\bm{\nu} \in I$, we define
  \begin{align*}
    \bm{\mu} * \bm{\nu} &= \sum_{M=(m_{ij}) \in \mathcal{H}_{\bm{\mu}, \bm{\nu}}} (m_{11}+m_{21}, \ldots, m_{1l}+m_{2l}) \in V \\
    \intertext{and}
    \bm{\mu} \,\overline{*}\, \bm{\nu} 
    &= \sum_{M=(m_{ij}) \in \mathcal{H}_{\bm{\mu}, \bm{\nu}}} (-1)^{l(\bm{\mu}) + l(\bm{\nu}) - l}\, (m_{11}+m_{21}, \ldots, m_{1l}+m_{2l}) \in V, 
  \end{align*}
  where $l$ is the number of columns of $M$.
  We extend $*$, $\overline{*} \colon I \times I \rightarrow V$ bilinearly on $V \times V$.
\end{definition}
For example, we have
\begin{displaymath}
 (1) * (2,3) = (1,2,3) + (2,1,3) + (2,3,1) + (3,3) + (2,4) 
\end{displaymath}
and
\begin{displaymath}
 (1) \,\overline{*}\, (2,3) = (1,2,3) + (2,1,3) + (2,3,1) - (3,3) - (2,4)
\end{displaymath}
from (\ref{eq2-25}).
\begin{proposition} 
 \label{prop2-40} 
  For any $v$, $w \in V$, we have 
  $d(v \, \overline{*} \, w) = d(v) * d(w)$.
\end{proposition}
\begin{proof}
 For a proof of this proposition, see~\cite[Proof of Proposition 2.4]{M}. 
\end{proof}
For $\bm{\mu} = (\mu_1,\ldots,\mu_p)$, $\bm{\nu} = (\nu_1,\ldots,\nu_q) \in I$, we define
\begin{align*}
  \bm{\mu} \, \# \, \bm{\nu} &= (\mu_1,\ldots,\mu_p,\nu_1,\ldots,\nu_q),\\
  \bm{\mu} \, \dot{\#} \, \bm{\nu} &= (\mu_1,\ldots,\mu_{p-1},\mu_p+\nu_1,\nu_2,\ldots,\nu_q)\\
  \intertext{and}
  \bm{\mu}^{+} &= (\mu_1,\ldots,\mu_{p-1},\mu_p +1).\\
\end{align*}
In addition, if $|\bm{\mu}| > 1$, we define
\begin{displaymath}
  \bm{\mu}^{-} =
  \begin{cases}
    (\mu_1,\ldots,\mu_{p-1},\mu_p-1) & \text{if $\mu_p > 1$}\\
    (\mu_1,\ldots,\mu_{p-1})         & \text{if $\mu_p = 1$}.
  \end{cases}
\end{displaymath}
We extend $\#$, $\dot{\#} \colon I \times I \rightarrow V$ bilinearly on
$V \times V$ and $\bmu \mapsto \bmu^{+}$ linearly on $V$. 
The following Propositions \ref{prop2-45} and \ref{prop2-50} are easily proved.
\begin{proposition} \label{prop2-45} 
 For any $\bmu$, $\bnu \in I$, we have
 $(\bmu \shp \bnu)^{*} = \bmu^{*} \dotshp \bnu^{*}$.
\end{proposition}
\begin{proposition} \label{prop2-50} 
 For any $\bmu$, $\bnu \in I$, we have\\[1mm]
 $(\mathrm{i})$ $u(\bmu \,\#\, \bnu) = u(\bmu) \,\#\, u(\bnu)$,\\[1mm]
 $(\mathrm{ii})$ 
 $u(\bmu \,\dot{\#}\, \bnu) = u(\bmu) \,\#\, u(\bnu) + u(\bmu) \,\dot{\#}\, u(\bnu)$,\\[1mm]
 $(\mathrm{iii})$ $d(\bmu \,\dot{\#}\, \bnu) = d(\bmu) \,\dot{\#}\, d(\bnu)$,\\[1mm]
 $(\mathrm{iv})$ 
 $d(\bmu \,\#\, \bnu) = d(\bmu) \,\#\, d(\bnu) + d(\bmu) \,\dot{\#}\, d(\bnu)$.
\end{proposition}
Here, we introduce the empty multi-index $\phi$.
We denote $\{\phi\} \cup I$ by $\tilde{I}$ and
$\Q \phi \oplus V$ by $\tilde{V}$.
It is natural to define
\begin{displaymath}
 \phi^{+} = (1) \quad \text{and} \quad (1)^{-} = \phi.
\end{displaymath}
We define $*(\phi)$, $\tau(\phi)$, $\sigma(\phi)$, $u(\phi)$ and
$d(\phi)$ to be $\phi$.
Moreover, we define for any $\bmu \in \tilde{I}$
\begin{displaymath}
 \phi * \bmu = \bmu * \phi = \bmu, \quad
 \phi \barast \bmu = \bmu \barast \phi = \bmu
\end{displaymath}
and
\begin{displaymath}
 \phi \shp \bmu = \bmu \shp \phi = \bmu, \quad
 \phi \dotshp \bmu = \bmu \dotshp \phi = \bmu.
\end{displaymath}
Proposition \ref{prop2-60} follows from
the definition of the harmonic product $*$.
\begin{proposition} \label{prop2-60}
 For any $\bmu = (\mu_1,\ldots,\mu_p)$, 
 $\bnu = (\nu_1,\ldots,\nu_q) \in I$, we have
 \begin{displaymath}
  \bmu * \bnu = (\bmu^{-} * \bnu) \,\#\, (\mu_p) + (\bmu * \bnu^{-}) \,\#\,
  (\nu_q) + (\bmu^{-} * \bnu^{-}) \,\#\, (\mu_p + \nu_q).
 \end{displaymath}
\end{proposition}
In the last of this section, we define bilinear mappings
$\circledast, \barcast \colon V \times V \to V$, which are used in the following
sections, and prove a proposition (Proposition \ref{prop2-70}). Lemma
\ref{lem2-70} is used in the proof of Proposition \ref{prop2-70}.
We note that
\begin{equation}
 d(\bmu) = \sum_{i=1}^p d(\mu_1,\ldots,\mu_{i-1}) \shp (\mu_i + \cdots +
  \mu_p) \label{eq2-60}
\end{equation}
for any $\bmu = (\muvec) \in I$, 
which is immediate from the definition of $d$.
\begin{lemma} 
 \label{lem2-70}
 For any $\bmu = (\muvec)$, $\bnu = (\nuvec) \in I$, we have
 \begin{align*}
  &d(\bmu) * d(\bnu) = \sum_{i=1}^p \Big(d(\mu_1,\ldots,\mu_{i-1}) *
  d(\bnu)\Big) \shp (\mu_i + \cdots + \mu_p)\\
  &+ \sum_{j=1}^q \Big(d(\bmu) * d(\nu_1,\ldots,\nu_{j-1})\Big) \shp
  (\nu_j + \cdots + \nu_q) \\
  &+ \sum_{i=1}^p \sum_{j=1}^q \Big(d(\mu_1,\ldots,\mu_{i-1}) *
  d(\nu_1,\ldots,\nu_{j-1})\Big)
  \shp (\mu_i + \cdots + \mu_p + \nu_j + \cdots + \nu_q).
 \end{align*}
\end{lemma}
\begin{proof}
 It follows from (\ref{eq2-60}) and Proposition \ref{prop2-60}.
\end{proof}
\begin{definition}
 For $\bmu = (\muvec)$, $\bnu = (\nuvec) \in I$, we define
 \begin{displaymath}
  \bmu \circledast \bnu = (\bmu^{-} * \bnu^{-}) \shp (\mu_p + \nu_q)
 \end{displaymath}
 and
 \begin{displaymath}
  \bmu \,\overline{\circledast}\, \bnu = (\bmu^{-} \,\overline{*}\, \bnu^{-}) \shp (\mu_p + \nu_q).
 \end{displaymath}
 We extend 
 $\circledast$, $\overline{\circledast} \colon I \times I \to V$ 
 bilinearly on $V \times V$.
\end{definition}
\begin{proposition} 
 \label{prop2-70}
 For any $v$, $w \in V$, we have
 $d(v \barcast w) = d(v) \circledast d(w)$.
\end{proposition}
\begin{proof}
 It suffices to show that the equality 
 $d(\bmu \barcast \bnu) = d(\bmu) \circledast d(\bnu)$
 holds for any $\bmu$, $\bnu \in I$. Let $\bmu = (\muvec)$, 
 $\bnu = (\nuvec) \in I$. Then we have
 \begin{equation}
  d(\bmu \barcast \bnu) 
  = \Bigl(d(\bmu^{-}) * d(\bnu^{-})\Bigr) \shp (\mu_p + \nu_q) 
  + \Bigl(d(\bmu^{-}) * d(\bnu^{-})\Bigr) \dotshp (\mu_p + \nu_q) \label{eq2-70}
 \end{equation}
 by the definition of $\barcast$, Propositions \ref{prop2-50} and \ref{prop2-40}. By Lemma \ref{lem2-70},
 we have
 \begin{align*}
  &\Bigl(d(\bmu^{-}) * d(\bnu^{-})\Bigr) \dotshp (\mu_p + \nu_q) \\
  &= \sum_{i=1}^{p-1} \Big(d(\mu_1,\ldots,\mu_{i-1}) * d(\bnu^{-})\Big)
  \shp (\mu_i + \cdots + \mu_p + \nu_q)\\
  &+ \sum_{j=1}^{q-1} \Big(d(\bmu^{-}) * d(\nu_1,\ldots,\nu_{j-1})\Big)
  \shp (\mu_p + \nu_j + \cdots + \nu_q)\\
  &+ \sum_{i=1}^{p-1} \sum_{j=1}^{q-1} \Big(d(\mu_1,\ldots,\mu_{i-1}) *
  d(\nu_1,\ldots,\nu_{j-1})\Big) \shp (\mu_i + \cdots + \mu_p + \nu_j +
  \cdots + \nu_q).
 \end{align*}
  Hence the right-hand side of (\ref{eq2-70}) equals 
  \begin{displaymath}
   \sum_{i=1}^p \sum_{j=1}^q \Big(d(\mu_1,\ldots,\mu_{i-1}) *
  d(\nu_1,\ldots,\nu_{j-1}) \Big) \shp (\mu_i + \cdots + \mu_p + \nu_j +
  \cdots + \nu_q),
  \end{displaymath}
  which is equal to $d(\bmu) \circledast d(\bnu)$ by (\ref{eq2-60}).
 \end{proof}
%
%
\section{A formula for finite multiple harmonic sums} \label{multiple-sum}
In this section, we prove a formula for finite multiple harmonic sums (Theorem \ref{th3-10}).
We note that $\mathbb{N}$ denotes the set of non-negative integers.
\begin{definition}
  (i) We define the difference $\Delta \colon \C^{\N} \to \C^{\N}$ by putting
  \begin{displaymath}
    (\Delta a)(n) = a(n) - a(n+1)
  \end{displaymath}
  for $a \in \C^{\N}$ and $n \in \mathbb{N}$.
  We denote the composition of $\Delta$ with itself $n$ times by $\Delta^n$.\\
  (ii) We define the inversion $\nabla \colon \C^{\N} \to \C^{\N}$ by putting
  \begin{displaymath}
    (\nabla a)(n) = (\Delta^n a)(0)
  \end{displaymath}
  for $a \in \C^{\N}$ and $n \in \mathbb{N}$.
\end{definition}
Let $\C[[x,y]]$ be the ring of formal power series over the complex
numbers in two indeterminates $x$ and $y$.
We put $\partial_x = \partial / \partial x$ and 
$\partial_y = \partial/\partial y$. For
\begin{displaymath}
 f(x,y) = \sum_{n,k=0}^{\infty} a(n,k) \, \frac{x^n y^k}{n! \, k!} \in \C[[x,y]],
\end{displaymath} 
we have
\begin{equation} \label{eq3-10}
 (\partial_x + \partial_y - 1) f(x,y) = \sum_{n,k=0}^{\infty} \{a(n+1,k) +
 a(n,k+1) - a(n,k)\} \, \frac{x^n y^k}{n! \, k!}.
\end{equation}
Therefore by the definition of the difference we have
\begin{equation} \label{eq3-20}
 (\partial_x + \partial_y - 1) \sum_{n,k=0}^{\infty} (\Delta^k a) (n) \,
 \frac{x^n y^k}{n! \, k!} = 0.
\end{equation}
\begin{lemma} \label{lem3-10}
 Let $f(x,y) \in \C[[x,y]]$. If $f(x,y)$ satisfies two conditions
 \begin{displaymath}
  (\partial_x + \partial_y - 1) f(x,y) = 0 \quad \text{and} \quad f(x,0)
  = 0,
 \end{displaymath}
 we have $f(x,y) = 0$.
\end{lemma}
\begin{proof}
 Let
 \begin{displaymath}
  f(x,y) = \sum_{n,k=0}^{\infty} a(n,k) \, \frac{x^n y^k}{n! \, k!}.
 \end{displaymath}
 By (\ref{eq3-10}), the conditions of the lemma are equivalent to
\begin{displaymath}
 \begin{cases}
  a(n+1,k) + a(n,k+1) = a(n,k) \quad (n,k \in \N), \\
  a(n,0) = 0 \quad (n \in \N).
 \end{cases}
\end{displaymath}
By induction on $k$, we obtain $a(n,k) = 0$ for any $n,k \in \N$.
\end{proof}
\begin{proposition} \label{prop3-10}
 Let $a \in \C^{\N}$. For any $n$, $k \in \N$, we have
\begin{displaymath}
 (\Delta^k(\nabla a))(n) = (\Delta^n a)(k).
\end{displaymath}
In particular, we obtain $\nabla(\nabla a) = a$ by setting $n = 0$.
\end{proposition}
\begin{proof}
 For a sequence $a \in \C^{\N}$, we put
\begin{displaymath}
 f_a(x,y) = \sum_{n,k=0}^{\infty} (\Delta^k a) (n) \,
 \frac{x^n y^k}{n! \, k!}.
\end{displaymath}
By (\ref{eq3-20}) we have 
$(\partial_x + \partial_y - 1) f_{\nabla a}(x,y) = 0$
and $(\partial_x + \partial_y - 1) f_{a}(y,x) = 0$.
We also have
\begin{displaymath}
 f_{\nabla a}(x,0) = \sum_{n=0}^{\infty}(\nabla a)(n) \frac{x^n}{n!} = f_a(0,x)
\end{displaymath}
by the definition of the inversion $\nabla$.
Therefore by Lemma \ref{lem3-10}, we obtain
$f_{\nabla a}(x,y) = f_a(y,x)$.
This completes the proof of Proposition \ref{prop3-10}.
\end{proof}
Let
\begin{displaymath}
 F(x) = \sum_{n=0}^{\infty}a(n) \frac{x^n}{n!}.
\end{displaymath}
By using Lemma \ref{lem3-10}, we can easily prove that
\begin{displaymath}
 F(x-y) e^y = \sum_{n,k=0}^{\infty} (\Delta^k a) (n) \,
 \frac{x^n y^k}{n! \, k!}.
\end{displaymath}
By comparing the coefficients of $x^n y^k$ on both sides, we obtain
the following proposition.
\begin{proposition} \label{prop3-20}
 For any $n$, $k \in \N$, we have
\begin{displaymath}
 (\Delta^k a)(n) = \sum_{i=0}^k (-1)^i \binom{k}{i} a(n+i).
\end{displaymath}
In particular, we obtain
\begin{displaymath}
 (\nabla a)(k) = \sum_{i=0}^k (-1)^i \binom{k}{i} a(i)
\end{displaymath}
for any $k \in \N$ by setting $n=0$.
\end{proposition}
Let $\bm{\mu} = (\mu_1,\ldots,\mu_p) \in I$ and $n \in \mathbb{N}$.
In this paper, we consider four kinds of finite multiple harmonic sums
\begin{align*}
  s_{\bm{\mu}}(n) &= \sum_{0 \le n_1 \le \cdots \le n_p = n}\frac{1}{(n_1+1)^{\mu_1} \cdots (n_p+1)^{\mu_p}},\\
  a_{\bm{\mu}}(n) &= \sum_{0 \le n_1 < \cdots < n_p = n}\frac{1}{(n_1+1)^{\mu_1} \cdots (n_p+1)^{\mu_p}},\\
  S_{\bm{\mu}}(n) &= \sum_{0 \le n_1 \le \cdots \le n_p < n}\frac{1}{(n_1+1)^{\mu_1} \cdots (n_p+1)^{\mu_p}}\\
\intertext{and}
  A_{\bm{\mu}}(n) &= \sum_{0 \le n_1 < \cdots < n_p < n}\frac{1}{(n_1+1)^{\mu_1} \cdots (n_p+1)^{\mu_p}}.
\end{align*}
We extend the mappings $I \ni \bm{\mu} \mapsto s_{\bm{\mu}}$,
$a_{\bm{\mu}}$, $S_{\bm{\mu}}$, 
$A_{\bm{\mu}} \in \mathbb{C}^{\mathbb{N}}$ linearly on $V$.
We introduced the empty multi-index $\phi$ in Section
\ref{multi-indices}. We define
\begin{displaymath}
 S_{\phi}(n) = A_{\phi}(n) = 1
\end{displaymath}
for any $n \in \N$.
\begin{proposition} 
\label{prop3-40} 
  \rm{(i)} For any $\bm{\mu} \in I$, we have
  $-(\Delta S_{\bm{\mu}}) = s_{\bm{\mu}}$ and $-(\Delta A_{\bm{\mu}}) = a_{\bm{\mu}}$.
  \rm{(ii)} For any $\bm{\mu} = (\muvec) \in I$ and any $n \in \mathbb{N}$, we have
  \begin{displaymath}
    s_{{\bm{\mu}}}(n) = S_{\mu_1,\ldots,\mu_{p-1}}(n+1)\frac{1}{(n+1)^{\mu_p}}
    \quad \text{and} \quad
    a_{{\bm{\mu}}}(n) = A_{\mu_1,\ldots,\mu_{p-1}}(n)\frac{1}{(n+1)^{\mu_p}},
  \end{displaymath}
  where we put $(\mu_1,\ldots,\mu_{p-1}) = \phi$ if $p = 1$.
\end{proposition}
\begin{proof}
 These are easily seen.
\end{proof}
\begin{proposition} \label{prop3-30} 
  \rm{(i)} For any $v \in V$, we have $s_v = a_{d(v)}$, $S_v = A_{d(v)}$, $a_v = s_{d^{-1}(v)}$ and $A_v = S_{d^{-1}(v)}$.\\
  \rm{(ii)} For any $v$, $w \in V$, we have 
  $A_v A_w = A_{v * w}$, $S_v S_w = S_{v \overline{*} w}$,
  $a_v a_w = a_{v \circledast w}$ and $s_v s_w = s_{v \barcast w}$.
\end{proposition}
\begin{proof}
  (i) It is easily seen that
  $s_{\bm{\mu}} = a_{d(\bm{\mu})}$ and $S_{\bm{\mu}} = A_{d(\bm{\mu})}$
  for $\bm{\mu} \in I$. The assertions of (i) follow immediately from these.\\
  (ii) It is easily seen that
  $A_{\bm{\mu}} A_{\bm{\nu}} = A_{\bm{\mu} * \bm{\nu}}$
  for any $\bm{\mu}$, $\bm{\nu} \in I$. Therefore the first assertion follows.
  In addition, for any $v$, $w \in V$ we have
  \begin{displaymath}
    S_v S_w = A_{d(v)} A_{d(w)} = A_{d(v)*d(w)} = S_{d^{-1}(d(v)*d(w))} = S_{v \overline{*} w}
  \end{displaymath}
  by Proposition \ref{prop2-40}.  
  The other assertions follow from Proposition \ref{prop3-40} (ii) and
  the definitions of $\circledast$ and $\barcast$.
\end{proof}
Let $m$ be a positive integer. 
For $\bm{\mu} = (\mu_1,\ldots,\mu_p)$, $\bm{\nu} = (\nu_1,\ldots,\nu_q) \in I_m$ and $n$, $k \in \mathbb{N}$, we also consider a finite sum
\begin{displaymath}
  s_{\bm{\mu},\bm{\nu}}(n,k) = \binom{n+k}{n}^{-1} \sum_{\substack{0 \le n_1 \le \cdots \le n_p = n\\ 0 \le k_1 \le \cdots \le k_q = k}}
  \frac{1}{(n_{i_1}+k_{j_1}+1) \cdots (n_{i_m}+k_{j_m}+1)},
\end{displaymath}
where
\begin{align*}
  (i_1,\ldots,i_m) &= (\underbrace{1,\ldots,1}_{\mu_1},\ldots,\underbrace{p,\ldots,p}_{\mu_p})\\
\intertext{and}
  (j_1,\ldots,j_m) &= (\underbrace{1,\ldots,1}_{\nu_1},\ldots,\underbrace{q,\ldots,q}_{\nu_q}).
\end{align*}
We note that
\begin{equation} 
\label{eq3-80}
  s_{\bm{\mu},\bm{\nu}}(n,0) = s_{\bm{\mu}}(n) \quad \text{and} \quad s_{\bm{\mu},\bm{\nu}}(0,k) = s_{\bm{\nu}}(k).
\end{equation}
\begin{lemma} \label{lem3-20}
  Let $\bm{\mu} = (\mu_1,\ldots,\mu_p)$, $\bm{\nu} = (\nu_1,\ldots,\nu_q) \in I$ and let $|\bm{\mu}|=|\bm{\nu}|=m$.
  Then for any $n$, $k \in \mathbb{N}$, we have
  \begin{displaymath}
      s_{\bm{\mu}^{-}, \bm{\nu}^{-}}(n,k) = 
      \begin{cases}
        (n+k+1) s_{\bm{\mu}, \bm{\nu}}(n,k) - k\, s_{\bm{\mu},
       \bm{\nu}}(n,k-1) & \!\!\! \text{if $\mu_p > 1$ and $\nu_q = 1$}\\
        (n+k+1) s_{\bm{\mu}, \bm{\nu}}(n,k) - n\, s_{\bm{\mu},
       \bm{\nu}}(n-1,k) & \!\!\! \text{if $\mu_p = 1$ and $\nu_q > 1$,}
      \end{cases}
  \end{displaymath}
  where we put $k\, s_{\bm{\mu}, \bm{\nu}}(n,k-1) = 0$ if $k=0$ and $n\, s_{\bm{\mu}, \bm{\nu}}(n-1,k)=0$ if $n=0$.
\end{lemma}
\begin{proof}
  We shall prove the lemma only for $\mu_p > 1$ and $\nu_q = 1$.
  We can prove the lemma for $\mu_p = 1$ and $\nu_q > 1$ similarly.
  We have
  \begin{align*}
    &(n+k+1) s_{\bm{\mu}, \bm{\nu}}(n,k) - k\, s_{\bm{\mu}, \bm{\nu}}(n,k-1)\\
    &= \binom{n+k}{n}^{-1} (n+k+1) 
       \sum_{\substack{0 \le n_1 \le \cdots \le n_p = n\\ 0 \le k_1 \le \cdots \le k_q = k}}
       \frac{1}{(n_{i_1} + k_{j_1} + 1) \cdots (n_{i_m} + k_{j_m} + 1)} \\
       & \quad - \binom{n+k}{n}^{-1} (n+k)
       \sum_{\substack{0 \le n_1 \le \cdots \le n_p = n\\ 0 \le k_1 \le \cdots \le k_q = k-1}}
       \frac{1}{(n_{i_1} + k_{j_1} + 1) \cdots (n_{i_m} + k_{j_m} + 1)}. \\
  \end{align*}
  Since $\mu_p > 1$ and $\nu_q = 1$, this equals
  \begin{align*}
    & \binom{n+k}{n}^{-1}
      \left\{ \sum_{\substack{0 \le n_1 \le \cdots \le n_p = n\\ 0 \le k_1 \le \cdots \le k_{q-1} \le k}}
      \frac{1}{(n_{i_1} + k_{j_1} + 1) \cdots (n_{i_{m-1}} + k_{j_{m-1}} + 1)} \right.\\
      & \left. \qquad\qquad - \sum_{\substack{0 \le n_1 \le \cdots \le n_p = n\\ 0 \le k_1 \le \cdots \le k_{q-1} \le k-1}}
      \frac{1}{(n_{i_1} + k_{j_1} + 1) \cdots (n_{i_{m-1}} + k_{j_{m-1}} + 1)} \right\}\\
    &= \binom{n+k}{n}^{-1}
       \sum_{\substack{0 \le n_1 \le \cdots \le n_p = n\\ 0 \le k_1 \le \cdots \le k_{q-1} = k}}
       \frac{1}{(n_{i_1} + k_{j_1} + 1) \cdots (n_{i_{m-1}} + k_{j_{m-1}} + 1)}.
  \end{align*}
  We have thus proved the lemma.
\end{proof}
\begin{theorem} 
\label{th3-10}
  For any $\bm{\mu} \in I$ and $n$, $k \in \mathbb{N}$ we have
  $(\Delta^k s_{\bm{\mu}})(n) = s_{\bm{\mu},\bm{\mu}^*}(n,k)$.
\end{theorem}
\begin{proof}
 For $\bm{\mu} \in I$, we put
 \begin{displaymath}
  f_{\bm{\mu}}(x,y) = \sum_{n,k=0}^{\infty}(\Delta^k s_{\bm{\mu}})(n)
  \frac{x^n y^k}{n! \, k!}
 \end{displaymath}
 and for $\bm{\mu}$, $\bm{\nu} \in I$ with $|\bm{\mu}| = |\bm{\nu}|$,
 we put
 \begin{displaymath}
  g_{\bm{\mu},\bm{\nu}}(x,y) = \sum_{n,k=0}^{\infty}
  s_{\bm{\mu},\bm{\nu}}(n,k) \frac{x^n y^k}{n! \, k!}.
 \end{displaymath}
 If we have
 \begin{equation}
  (\partial_x + \partial_y - 1) g_{\bmu,\bmu^{*}}(x,y) = 0, \label{eq3-90}
 \end{equation}
 then by Lemma \ref{lem3-10} and (\ref{eq3-80}) we have
 $f_{\bmu}(x,y) = g_{\bmu,\bmu^{*}}(x,y)$.
 Therefore we have only to prove (\ref{eq3-90}). If we put
 \begin{displaymath}
  \xi = x \partial_x + y \partial_y + 1 - y \quad \text{and} \quad \eta
  = x \partial_x + y \partial_y + 1 - x,
 \end{displaymath}
 by Lemma \ref{lem3-20} we have
 \begin{displaymath}
    g_{\mu^{-},\nu^{-}} = 
    \begin{cases}
      \xi  g_{\mu,\nu} & \text{if $\mu_p > 1$ and $\nu_q = 1$}\\
      \eta g_{\mu,\nu} & \text{if $\mu_p = 1$ and $\nu_q > 1$}.
    \end{cases}
  \end{displaymath}
 We note that
 \begin{displaymath}
  (\partial_x + \partial_y - 1)\xi - \xi (\partial_x + \partial_y + 1) = \partial_x + \partial_y - 1
 \end{displaymath}
 and
 \begin{displaymath}
  (\partial_x + \partial_y - 1)\eta - \eta (\partial_x + \partial_y + 1) = \partial_x + \partial_y - 1.
 \end{displaymath}
 We prove (\ref{eq3-90}) by induction on $|\bmu|$. For $|\bmu| = 1$
 (i.e. $\bmu = (1)$), we have
 \begin{displaymath}
  g_{\bmu,\bmu^{*}}(x,y) = \sum_{n,k=0}^{\infty} \frac{x^n y^k}{(n+k+1)!}.
 \end{displaymath}
 Therefore we can easily check that (\ref{eq3-90}) holds in this
 case. We assume that the claim (\ref{eq3-90}) is true for 
 $|\bmu| = m$. Let $|\bmu| = m+1$, $\bmu = (\mu_1,\ldots,\mu_p)$ and
 $\bmu^{*} = (\mu^{*}_1,\ldots,\mu^{*}_q)$.
 If $\mu_p > 1$, then we have $\mu^{*}_q = 1$ by (\ref{eq2-10}).
 Therefore we obtain
 \begin{align*}
  &(\xi + 1)(\partial_x + \partial_y - 1) g_{\bmu,\bmu^{*}} =
  (\partial_x + \partial_y - 1) \xi g_{\bmu,\bmu^{*}}
  = (\partial_x + \partial_y - 1) g_{\bmu^{-},(\bmu^{*})^{-}}\\
  &= (\partial_x + \partial_y - 1) g_{\bmu^{-},(\bmu^{-})^{*}} = 0
 \end{align*}
 by the hypothesis of the induction.
 If $\mu_p = 1$, then we obtain
 $(\eta + 1)(\partial_x + \partial_y - 1) g_{\bmu,\bmu^{*}} = 0$
 by a similar argument.
 We can easily check that mappings 
 $\xi + 1$, $\eta + 1 \colon \C[[x,y]] \to \C[[x,y]]$ are injections.
 Thus we have proved (\ref{eq3-90}).
\end{proof}
\begin{corollary} 
 \label{cor3-10}
 For any $\bmu \in I$, we have $\nabla s_{\bmu} = s_{\bmu^*}$.
\end{corollary}
\begin{proof}
  If we set $n=0$ in Theorem \ref{th3-10}, we obtain the result by (\ref{eq3-80}). 
\end{proof}
We define $T \colon \C^{\N} \to \C^{\N}$ by
\begin{displaymath}
 (Ta)(n) = a(n+1) \quad (a \in \C^{\N},\, n \in \N),
\end{displaymath}
and then we have $\Delta = \nabla T \nabla$.
\begin{corollary} 
 \label{cor3-20} 
 For any $\bmu \in I$ and $n \in \mathbb{N}$, we have 
 \begin{displaymath}
  (\nabla S_{\bmu})(n) =
  \begin{cases}
   0 & \text{if $n=0$}\\
   -s_{\bmu^{*}}(n-1) & \text{if $n \ge 1$}.
  \end{cases}
 \end{displaymath}
\end{corollary}
\begin{proof}
 By Proposition \ref{prop3-40} (i) and Corollary \ref{cor3-10} we have
 $T \nabla S_{\bmu} = - s_{\bmu^{*}}$,
 which completes the proof.
\end{proof}
The following estimation is used in Section \ref{multiple-zeta}.
\begin{proposition} 
 \label{prop3-50} 
 Let $k \in \N$ and let $\bmu = (\mu_1,\ldots,\mu_p) \in I$.
 For any $\varepsilon > 0$, we have
 \begin{displaymath}
  (\Delta^k s_{\bmu})(n) =
  O\left(\frac{1}{n^{k+\mu_p-\varepsilon}}\right) \quad (n \to \infty).
 \end{displaymath} 
\end{proposition}
\begin{proof}
 Let $|\bmu| = m$ and let $\bmu^{*} = (\mu^{*}_1,\ldots,\mu^{*}_q)$. By
 Theorem \ref{th3-10} we have
 \begin{displaymath}
  (\Delta^k s_{\bmu})(n) = \binom{n+k}{n}^{-1} \sum_{\substack{0 \le n_1
  \le \cdots \le n_p = n\\ 0 \le k_1 \le \cdots \le k_q = k}} \frac{1}{(n_{i_1}+k_{j_1}+1)\cdots(n_{i_m}+k_{j_m}+1)},
 \end{displaymath}
 where
 \begin{displaymath}
  (i_1,\ldots,i_m) =
  (\underbrace{1,\ldots,1}_{\mu_1},\ldots,\underbrace{p,\ldots,p}_{\mu_p}), \quad
  (j_1,\ldots,j_m) =
  (\underbrace{1,\ldots,1}_{\mu^{*}_1},\ldots,\underbrace{q,\ldots,q}_{\mu^{*}_q}).
 \end{displaymath}
 The number of elements of the set
 \begin{displaymath}
  \{(k_1,\ldots,k_{q-1}) \in \N^{q-1} \,|\, 0 \le k_1 \le \cdots \le
  k_{q-1} \le k\}
 \end{displaymath} 
 is bounded by $(k+1)^{q-1}$. Hence we have
 \begin{align*}
  &\sum_{\substack{0 \le n_1
  \le \cdots \le n_p = n\\ 0 \le k_1 \le \cdots \le k_q = k}}
  \frac{1}{(n_{i_1}+k_{j_1}+1)\cdots(n_{i_m}+k_{j_m}+1)} \notag \\
  &\le \sum_{\substack{0 \le n_1 \le \cdots \le n_p = n\\ 0 \le k_1 \le
  \cdots \le k_q = k}}
  \frac{1}{(n_{1}+1)^{\mu_1}\cdots(n_{p}+1)^{\mu_p}} \\
  &\le (k+1)^{q-1} \sum_{0 \le n_1 \le \cdots \le n_p =
  n} \frac{1}{(n_1+1)^{\mu_1} \cdots (n_p+1)^{\mu_p}}\\
  &\le \frac{(k+1)^{q-1}}{(n+1)^{\mu_p}} \left(\sum_{0 \le n_1 \le
  n}\frac{1}{n_1+1}\right) \cdots \left(\sum_{0 \le n_{p-1} \le
  n}\frac{1}{n_{p-1}+1}\right)\\
  &= O \left(\frac{\log^{p-1}n}{n^{\mu_p}}\right) \quad (n \to \infty),
 \end{align*}
 which completes the proof.
\end{proof}
Moreover, we use the inequality
\begin{equation}
 s_{\mu_1,\ldots,\mu_p}(n) 
 \ge \frac{1}{(n+1)^{\mu_p}} \quad (n \in \N)
 \label{eq3-140} 
\end{equation}
in Section \ref{multiple-zeta}.
%
%
\section{The Newton series} \label{The-Newton-series}
In this section, we state general properties of the Newton series needed
later. We do not
prove Propositions \ref{prop4-5}, \ref{prop4-10} and \ref{prop4-20}. For
proofs of these propositions, see~\cite{G,I,Mi}.
For a sequence $a \in \C^{\N}$, the series with complex variable $z$
\begin{equation}
 f(z) = \sum_{n=0}^{\infty}(-1)^n (\nabla a)(n)\binom{z}{n},
 \quad \binom{z}{n} = \frac{z(z-1)\cdots (z-n+1)}{n!} \label{eq4-1}
\end{equation}
is called the Newton series which interpolates a sequence 
$a$. In fact, for $z = n \in \N$, the sum of this series is
\begin{displaymath}
 \sum_{k=0}^{n}(-1)^k (\nabla a)(k)\binom{n}{k} = (\nabla^2 a)(n)
 = a(n)
\end{displaymath}
by Propositions \ref{prop3-20} and \ref{prop3-10}. 
We first state known results for the convergence of the Newton series.
\begin{proposition}
 \label{prop4-5}
 Let $a \in \C^{\N}$ and $z_0 \in \C \setminus \N$. If the Newton series
 \begin{displaymath}
  \sum_{n=0}^{\infty} (-1)^n a(n) \binom{z}{n}
 \end{displaymath}
 converges at $z = z_0$, then the series converges uniformly in the
 wider sense in the region $\re z > \re z_0$.
\end{proposition}
By Proposition \ref{prop4-5}, there exists the unique 
$\rho \in \R \cup \{\pm \infty\}$ such that the series $f(z)$ defined by 
(\ref{eq4-1}) converges if $\re z > \rho$ and does not converge if
$\re z < \rho$ and $z \notin \N$. We call $\rho$ the abscissa of
convergence of the series $f(z)$. The function $f(z)$ is analytic in the region $\re z > \rho$.
\begin{proposition} 
 \label{prop4-10}
 Let $a \in \C^{\N}$ and $z \in \C \setminus \N$. Then the Newton series
 \begin{displaymath}
  \sum_{n=0}^{\infty} (-1)^n a(n) \binom{z}{n} 
 \end{displaymath}
 converges $($resp. converges absolutely$)$ if and only if the Dirichlet series
 \begin{displaymath}
  \sum_{n=0}^{\infty} \frac{a(n)}{(n+1)^{z+1}} 
 \end{displaymath}
 converges $($resp. converges absolutely$)$.
\end{proposition}
The following corollary is immediate from Proposition \ref{prop4-10}.
\begin{corollary} 
\label{cor4-10}
 Let $a \in \C^{\N}$ and $\varepsilon \in \R$. If we have
 $a(n) = O(n^{\varepsilon})$ $(n \to \infty)$,
 then the Newton series
 \begin{displaymath}
  \sum_{n=0}^{\infty}(-1)^n a(n)\binom{z}{n}
 \end{displaymath}
 converges absolutely for any $z \in \C$ with $\re z > \varepsilon$.
\end{corollary}
The function $f(z)$ defined by a Newton series is uniquely determined by
$f(n)$ for all sufficiently large integers $n$.
\begin{proposition} 
 \label{prop4-20} 
 Let $a \in \C^{\N}$ and let the Newton series
 \begin{displaymath}
  f(z) = \sum_{n=0}^{\infty} (-1)^n a(n) \binom{z}{n}
 \end{displaymath}
 have the abscissa of convergence $\rho$. If there exists $N \in \N$
 such that $f(n)=0$ for any integers $n \ge N$,
 then we have $f(z) = 0$ for any $\re z > \rho$.
\end{proposition}
The product of two Newton series does not necessarily have a Newton
series expansion.
We give a sufficient condition for the product of two Newton series to
be expanded to a Newton series. We need a lemma.
\begin{lemma} 
 \label{lem4-50}
 Let $a \in \C^{\N}$ and let the Newton series
\begin{displaymath}
 f(z) = \sum_{n=0}^{\infty} (-1)^n a(n)\binom{z}{n}
\end{displaymath}
have the abscissa of convergence $\rho$. Let $k \in \N$. Then for any $\re z
 > \rho + k$, we have
\begin{displaymath}
 (-1)^k \binom{z}{k} f(z) =
 \sum_{n=k}^{\infty} (-1)^n \binom{n}{k}(\Delta^k a)(n-k)\binom{z}{n}.
\end{displaymath}
\end{lemma}
\begin{proof}
 For any $z \in \C$ and $n$, $k \in \N$, we have
 \begin{displaymath}
  \binom{z}{n} = \sum_{i=0}^{k} \binom{k}{i} \binom{z-k}{n-i}.
 \end{displaymath}
 Hence we have for any $\re z > \rho$
 \begin{align*}
  &(-1)^k \binom{z}{k} f(z) \\
  &= (-1)^k \sum_{i=0}^k \binom{k}{i} \sum_{n=i}^{\infty}(-1)^n a(n)
  \binom{z}{k} \binom{z-k}{n-i} \notag \\
  &= (-1)^k \sum_{i=0}^k \binom{k}{i} \sum_{n=i}^{\infty}(-1)^n a(n)
  \binom{k+n-i}{k}\binom{z}{k+n-i} \notag \\
  &= (-1)^k \sum_{i=0}^k \binom{k}{i} \sum_{n=0}^{\infty}(-1)^{n+i} a(n+i)
  \binom{k+n}{k}\binom{z}{k+n}. 
 \end{align*}
 By Proposition \ref{prop3-20}, this is equal to
 \begin{displaymath}
  (-1)^k \sum_{n=0}^{\infty}(-1)^n (\Delta^k a)(n)
  \binom{k+n}{k}\binom{z}{k+n},
 \end{displaymath}
 which completes the proof.
\end{proof}
\begin{proposition} 
 \label{prop4-30} 
 Let $a$, $b \in \C^{\N}$. Let the Newton series
 \begin{displaymath}
  f(z)=\sum_{n=0}^{\infty} (-1)^n a(n)\binom{z}{n} \quad \text{and} \quad
  g(z)=\sum_{n=0}^{\infty} (-1)^n b(n) \binom{z}{n}
 \end{displaymath}
 have the abscissas of convergence $\rho_a$ and $\rho_b$,
 respectively. 
 Let $\varepsilon > 0$.
 We assume that sequences $a$ and $b$ satisfy the
 following conditions$:$\\
 $(\mathrm{i})$ The values $a(n)$ and $(\Delta^k b)(n)$ are non-negative for any
 $n$, $k \in \N$.\\ 
 $(\mathrm{ii})$ $\rho_a < 0$.\\
 $(\mathrm{iii})$ For any $k \in \N$ we have
 $(\Delta^k b)(n) = O(n^{-k-\varepsilon})$ 
 $(n \to \infty)$. \\
 Then the product $f(z)g(z)$
 is expressed by a Newton series in the half-plane $\re z > \max\{\rho_a, -\varepsilon\}$.
\end{proposition}
\begin{proof}
 Let $k \in \N$. By Lemma \ref{lem4-50}, it holds that 
 \begin{equation}
  (-1)^k \binom{z}{k} g(z) =
   \sum_{n=k}^{\infty}(-1)^n\binom{n}{k}(\Delta^k b)(n-k)\binom{z}{n} \label{eq4-80}
 \end{equation}
  for any $\re z > \rho_b + k$.
  By setting $k=0$ in the condition (iii), we have
  $b(n) = O(n^{-\varepsilon})$ $(n \to \infty)$.
  According to Corollary \ref{cor4-10}, the Newton series $g(z)$
 converges in $\re z > -\varepsilon$.
  Again by the condition (iii), we have
 \begin{displaymath}
  \binom{n}{k} (\Delta^k b)(n-k) =
  O(n^{-\varepsilon}) \quad (n \to \infty).
 \end{displaymath}
 So the right-hand side of (\ref{eq4-80}) converges in 
 $\re z > -\varepsilon$. 
 Thus we see that the equality (\ref{eq4-80}) is
 valid for any $\re z > -\varepsilon$. Let $\re z > \max\{\rho_a,
 -\varepsilon\}$.
 We multiply the both sides of (\ref{eq4-80}) by $a(k)$ and take the sum
 over all values of $k$ from $0$ to $\infty$. Consequently, we obtain
 \begin{equation}
  f(z)g(z)
  = \sum_{k=0}^{\infty} \sum_{n=k}^{\infty} a(k)(-1)^n
  \binom{n}{k}(\Delta^k b)(n-k)\binom{z}{n}. \label{eq4-90}
 \end{equation}
 If we take $z$ as a negative real number with $z > \max\{\rho_a,
 -\varepsilon\}$, then each term of the right-hand side of (\ref{eq4-90}) is
 non-negative from the condition (i). 
 Hence the right-hand side of (\ref{eq4-90}) is equal to 
 \begin{displaymath}
  \sum_{n=0}^{\infty}(-1)^n \left\{\sum_{k=0}^n
  a(k)\binom{n}{k}(\Delta^k b)(n-k)\right\}\binom{z}{n}.
 \end{displaymath}
 Thus we have completed the proof.
\end{proof}
The Taylor expansion of a Newton series is given by the following
proposition. We note that the definition of the sequence
$a_{\bmu}$ $(\bmu \in I)$ is given in Section \ref{multiple-sum}.
\begin{proposition} 
 \label{prop4-40} 
 Let $a \in \C^{\N}$ be any sequence. We suppose that the abscissa of convergence of the
 Newton series
 \begin{displaymath}
  f(z) = \sum_{n=0}^{\infty}(-1)^n a(n)\binom{z}{n}
 \end{displaymath}
 is negative. Then we have
 \begin{displaymath}
  f(z) = a(0) + \sum_{m=1}^{\infty}(-1)^m \left\{\sum_{n=1}^{\infty}a(n)a_{\underbrace{1,\ldots,1}_m}(n-1)\right\}z^m
 \end{displaymath}
 in some neighborhood of $0$.
\end{proposition}
\begin{proof}
 Let $m$, $n \ge 1$. Then we have
 \begin{align*}
  \frac{1}{m!} \frac{d^m}{dz^m} \binom{z}{n} \biggm|_{z=0} 
  &= \binom{z}{n} \sum_{0 \le n_1 < \cdots < n_m < n}
  \frac{1}{(z-n_1)\cdots(z-n_m)} \biggm|_{z=0} \\
  &= (-1)^{n+m}\, a_{\underbrace{1,\ldots,1}_m}(n-1),
 \end{align*}
 which completes the proof.
\end{proof}
%
%
\section{A class of relations among multiple zeta values} \label{multiple-zeta}
In this section, we consider the Newton series which interpolates the
sequence $S_{\bmu}$ and derive a class of relations among MZV's. 
For $\bmu \in I$, we define
\begin{displaymath}
 f_{\bmu}(z) =
 \sum_{n=0}^{\infty}(-1)^n (\nabla s_{\bmu})(n) \binom{z}{n}
 \quad \text{and} \quad
 F_{\bmu}(z) = 
 \sum_{n=0}^{\infty}(-1)^n (\nabla S_{\bmu})(n) \binom{z}{n}.
\end{displaymath}
The abscissas of convergence of $f_{\bmu}(z)$ and $F_{\bmu}(z)$ are given by Proposition
\ref{prop5-10}. 
\begin{proposition} 
\label{prop5-10}
 Let $\bmu \in I$ and let
 \begin{displaymath}
  \bmu = (\mu_1,\ldots,\mu_r,\underbrace{1,\ldots,1}_l) \quad (\mu_r \ge
  2,\, r \ge 0,\, l \ge 0).
 \end{displaymath}
 Then the abscissas of convergence of the Newton series $f_{\bmu}(z)$
 and $F_{\bmu}(z)$ are both $-l-1$.
\end{proposition}
\begin{proof}
 We prove only the case $F_{\bmu}(z)$. The case $f_{\bmu}(z)$ is proved similarly.
 Let $\bmu^{*} = (\mu^{*}_1,\ldots,\mu^{*}_q)$. Then we have $\mu^{*}_q = l+1$ by
 (\ref{eq2-10}). By Corollary \ref{cor3-20} and 
 Proposition \ref{prop3-50}, we obtain
 \begin{displaymath}
  (\nabla S_{\bmu})(n) = O\left(\frac{1}{n^{l+1-\varepsilon}}\right)
  \quad (n \to \infty)
 \end{displaymath}
 for any $\varepsilon > 0$. Therefore $F_{\bmu}(z)$
 converges for any $\re z > -l-1$ by Corollary \ref{cor4-10}. 
 Since we have
 \begin{displaymath}
  s_{\bmu^{*}}(n) \ge \frac{1}{(n+1)^{l+1}} \quad (n \in N)
 \end{displaymath}
 by (\ref{eq3-140}), the Dirichlet series
 \begin{displaymath}
  \sum_{n=0}^{\infty} \frac{(\nabla S_{\bmu})(n)}{(n+1)^{z+1}} = -
  \sum_{n=1}^{\infty} \frac{s_{\bmu^{*}}(n-1)}{(n+1)^{z+1}}
 \end{displaymath}
 diverges at $z=-l-1$. Hence the assertion follows from Proposition \ref{prop4-10}.
\end{proof}
More generally, we define
\begin{displaymath}
 F_v(z) = \sum_{n=0}^{\infty} (-1)^n (\nabla S_{v})(n) \binom{z}{n}
\end{displaymath}
for any $v \in V$.
By Proposition \ref{prop5-10}, $F_v(z)$ converges at least in the region 
$\re z > -1$ and therefore it can be expanded in a Taylor series
at $z = 0$. We shall prove that the Taylor coefficients of $F_v(z)$ is
expressed by MZV's. 
We denote the $\Q$-vector
space whose basis is $I^{+} := \{\bmu^{+} \,|\, \bmu \in I\}$ by
$V^{+}$. We define $\Q$-linear mappings $\zeta$ and 
$\overline{\zeta}$ from $V^{+}$ to $\R$ by
\begin{displaymath}
 \zeta(v) = \sum_{n=0}^{\infty} a_v(n) \quad (v \in V^{+}) \quad \text{and} \quad
 \overline{\zeta}(v) = \sum_{n=0}^{\infty} s_v(n) \quad (v \in V^{+}),
\end{displaymath}
respectively. It is easily seen that
$\overline{\zeta}(v) = \zeta(d(v))$ and
$\zeta(v) = \overline{\zeta}(d^{-1}(v))$
from Proposition \ref{prop3-30} (i).
\begin{proposition} 
 \label{prop5-20} 
 Let $v \in V$. Then we have
 \begin{displaymath}
  F_v(z) = \sum_{m=1}^{\infty} (-1)^{m-1} \, \zeta \big(d\!*\!(v) \circledast
  (\underbrace{1,\ldots,1}_m) \big) z^m
 \end{displaymath}
 in some neighborhood of $0$.
\end{proposition}
\begin{proof}
 By Proposition \ref{prop4-40} and Corollary \ref{cor3-20}, we have
 \begin{displaymath}
  F_{v}(z)
  = \sum_{m=1}^{\infty} (-1)^{m-1} \left\{\sum_{n=1}^{\infty}
  s_{v^{*}}(n-1) \, a_{\underbrace{1,\ldots,1}_m}(n-1)\right\} z^m.
 \end{displaymath}
 Since
 \begin{displaymath}
  \sum_{n=0}^{\infty} s_{v^{*}}(n) \,
  a_{\underbrace{1,\ldots,1}_m}(n) 
  = \sum_{n=0}^{\infty} a_{d*(v) \, \circledast \,
  (\underbrace{1,\ldots,1}_m)}(n)\\
  = \zeta \big(d\!*\!(v) \circledast (\underbrace{1,\ldots,1}_m) \big)
 \end{displaymath}
 by Proposition \ref{prop3-30}, the proposition is proved.
\end{proof}
In Section \ref{multiple-sum}, we stated the harmonic product relations
for the finite multiple harmonic sums $S_v(n)$ (Proposition \ref{prop3-30} (ii)).
The functions $F_v(z)$ also satisfy the harmonic product relations which come
from those for $S_v(n)$.
\begin{theorem} 
\label{th5-10} 
 Let $v$, $w \in V$. Then we have
 \begin{equation}
  F_v(z) F_w(z) = F_{v \barast w}(z) \label{eq5-10} 
 \end{equation}
 for any $z \in \C$ with $\re z > -1$.
\end{theorem}
\begin{proof}
 Let $\bmu$, $\bnu \in I$.
 By Lemma \ref{lem4-50}, the function $(z+1) f_{\bmu \shp (1)}(z)$ is
 expanded to a Newton series. Therefore we have
 \begin{displaymath}
  (z+1)f_{\bmu \shp (1)}(z) = F_{\bmu}(z+1) \quad (\re z > -2)
 \end{displaymath}
 by Propositions \ref{prop3-40} (ii) and \ref{prop4-20}. Since
 $f_{\bmu \shp (1)}(z) f_{\bnu \shp (1)}(z)$ is expanded to a Newton
 series by Propositons \ref{prop3-50} and \ref{prop4-30}, the product
 $F_{\bmu}(z+1) F_{\bnu}(z+1)$ is also expanded to a Newton series. Hence we
 have
 \begin{displaymath}
  F_{\bmu}(z+1) F_{\bnu}(z+1) = F_{\bmu * \bnu}(z+1) \quad (\re z > -2),
 \end{displaymath}
 which completes the proof.
\end{proof}
Now we are in a stage to get the relations among MZV's.
\begin{corollary} 
\label{cor5-10} 
 For any $v$, $w \in V$ and any integer $m \ge 1$, we have
 \begin{displaymath}
  \sum_{\substack{k+l=m\\ k,l \ge 1}} \zeta \big(u \sigma (v) \circledast
  (\underbrace{1,\ldots,1}_k)\big) \, \zeta \big(u \sigma (w) \circledast
  (\underbrace{1,\ldots,1}_l)\big) = \zeta \big(u \sigma (v * w) \circledast (\underbrace{1,\ldots,1}_m)\big).
 \end{displaymath}
\end{corollary}
\begin{proof}
 We compare the Taylor coefficients in (\ref{eq5-10}). By Proposition
 \ref{prop5-20}, we have
 \begin{multline*}
  \sum_{\substack{k+l=m\\ k,l \ge 1}} \zeta \big(d\!*\!(v)
  \circledast (\underbrace{1,\ldots,1}_k)\big) \,
  \zeta \big(d\!*\!(w) \circledast (\underbrace{1,\ldots,1}_l)\big)\\
  = - \zeta \big(d\!*\! (v \barast w) \circledast (\underbrace{1,\ldots,1}_m)\big).
 \end{multline*}
 If we replace $v$ and $w$ by $d^{-1}(v)$ and $d^{-1}(w)$, respectively,
 we obtain the result since we have 
 $d^{-1}(v) \barast d^{-1}(w) = d^{-1} (v * w)$ and $d * d^{-1} = - u \sigma$
 by Propositions \ref{prop2-40} and \ref{prop2-20}.
\end{proof}
We define
\begin{displaymath}
 \zeta^{+}(v) = \zeta(v^{+}) \quad \text{and} \quad
 \overline{\zeta}{}^{+}(v) = \overline{\zeta}(v^{+})
\end{displaymath}
for any $v \in V$. Then it is easily seen that
\begin{equation}
 \overline{\zeta}{}^{+}(v) = \zeta^{+}\big(d(v)\big). 
 \label{eq5-20} 
\end{equation}
We denote the $\Q$-vector space generated by $v*w$ for all 
$v$, $w \in V$ by $V*V$.
The following corollary follows from Corollary \ref{cor5-10} on setting $m=1$.
\begin{corollary} 
 \label{cor5-20} 
 We have $u\sigma (V * V) \subset \ker \zeta^{+}$.
\end{corollary}
\begin{remark} 
 \label{rem5-10} 
 We can state Corollary \ref{cor5-20} as a claim for
 $\overline{\zeta}{}^{+}$. By (\ref{eq5-20}) we have
 $\overline{\zeta}{}^{+}\big(d^{-1}u\sigma (v*w)\big) = 0$
 for any $v$, $w \in V$.
 If we replace $v$ by $d(v)$ and $w$ by $d(w)$, then we obtain
 \begin{displaymath}
  \overline{\zeta}{}^{+}\big( (v \barast w)^{*}\big) = 0
 \end{displaymath}
 from Propositions \ref{prop2-40} and \ref{prop2-20}.
\end{remark}
\begin{remark}
 This remark was suggested by the referee. By a calculation similar to
 \begin{displaymath}
  A_2(n) 
  = \sum_{0 \le n_1 < n} \int_0^1 \! \int_0^1 (x_1 x_2)^{n_1} \, dx_1
  dx_2
  = \int_0^1 \! \int_0^1 \frac{1-(x_1 x_2)^n}{1-x_1 x_2} \, dx_1 dx_2,
 \end{displaymath}
 for any $\bmu = (\muvec) \in I$ and $n \in \N$ we obtain
 \begin{multline*}
  A_{\bmu}(n) = \int_0^1 \! \cdots \! \int_0^1
  \sum_{\varepsilon_1,\ldots,\varepsilon_{p-1} = 0,1}
  (-1)^{\varepsilon_1 + \cdots + \varepsilon_{p-1}} \\
  \times \frac{1- (y_1^{\varepsilon_1 \cdots \varepsilon_{p-1}}
  y_2^{\varepsilon_2 \cdots \varepsilon_{p-1}} \cdots
  y_{p-1}^{\varepsilon_{p-1}} y_p)^n}{(1-y_1)(1-y_1^{\varepsilon_1} y_2)
  \cdots
  (1-y_1^{\varepsilon_1 \cdots \varepsilon_{p-1}} y_2^{\varepsilon_2
  \cdots \varepsilon_{p-1}} \cdots y_{p-1}^{\varepsilon_{p-1}} y_p)} 
  \, dx_1 \cdots dx_{|\bmu|},
 \end{multline*}
 where 
 $y_i = x_{\mu_1+\cdots+\mu_{i-1}+1} \cdots x_{\mu_1+\cdots+\mu_{i-1}+\mu_i}$ $(1 \le i \le p)$.
 Therefore we conjecture that
 \begin{multline*}
  G_{\bmu}(z) := \sum_{n=0}^{\infty}(-1)^n (\nabla A_{\bmu})(n)
  \binom{z}{n} 
  = \int_0^1 \! \cdots \! \int_0^1
  \sum_{\varepsilon_1,\ldots,\varepsilon_{p-1} = 0,1}
  (-1)^{\varepsilon_1 + \cdots + \varepsilon_{p-1}} \\
  \times \frac{1- (y_1^{\varepsilon_1 \cdots \varepsilon_{p-1}}
  y_2^{\varepsilon_2 \cdots \varepsilon_{p-1}} \cdots
  y_{p-1}^{\varepsilon_{p-1}} y_p)^z}{(1-y_1)(1-y_1^{\varepsilon_1} y_2)
  \cdots
  (1-y_1^{\varepsilon_1 \cdots \varepsilon_{p-1}} y_2^{\varepsilon_2
  \cdots \varepsilon_{p-1}} \cdots y_{p-1}^{\varepsilon_{p-1}} y_p)}
  \, dx_1 \cdots dx_{|\bmu|}
 \end{multline*}
 for $\re z > 0$. Since
 \begin{displaymath}
  G_{\bmu}(z) = \sum_{m=1}^{\infty}(-1)^m \zeta \bigl(u\sigma(\bmu)
  \circledast (\underbrace{1,\ldots,1}_m) \bigr) z^m,
 \end{displaymath}
 we obtain
 \begin{multline*}
  \zeta^{+} \bigl(u\sigma(\bmu) \bigr) = \int_0^1 \! \cdots \! \int_0^1
  \sum_{\varepsilon_1,\ldots,\varepsilon_{p-1} = 0,1}
  (-1)^{\varepsilon_1 + \cdots + \varepsilon_{p-1}} \\
  \times \frac{\log(y_1^{\varepsilon_1 \cdots \varepsilon_{p-1}}
  y_2^{\varepsilon_2 \cdots \varepsilon_{p-1}} \cdots
  y_{p-1}^{\varepsilon_{p-1}} y_p)}{(1-y_1)(1-y_1^{\varepsilon_1} y_2)
  \cdots
  (1-y_1^{\varepsilon_1 \cdots \varepsilon_{p-1}} y_2^{\varepsilon_2
  \cdots \varepsilon_{p-1}} \cdots y_{p-1}^{\varepsilon_{p-1}} y_p)} 
  \, dx_1 \cdots dx_{|\bmu|}.
 \end{multline*}
 It seems to be possible to prove Corollary \ref{cor5-20} by standard
 manipulation on integrals. For example, if $\bmu = \bnu = (1)$, we have
 \begin{align*}
  &\zeta^{+} \bigl(u\sigma(\bmu * \bnu) \bigr) 
  = \zeta^{+} \bigl(u\sigma(2) \bigr) + 2 \zeta^{+} \bigl(u\sigma(1,1) \bigr) \\
  &= \int_0^1 \! \int_0^1 \left\{\frac{\log(x_1x_2)}{1-x_1x_2} +
  \frac{2\log x_2}{(1-x_1)(1-x_2)} -
  \frac{2\log(x_1x_2)}{(1-x_1)(1-x_1x_2)}\right\} \, dx_1 dx_2.
 \end{align*}
 Since the integrand is equal to
 \begin{displaymath}
  - \frac{1+x_1}{(1-x_1)(1-x_1x_2)} \log x_1 + \frac{1+x_2}{(1-x_2)(1-x_1x_2)} \log x_2,
 \end{displaymath}
 the right-hand side of the above equation is zero. This remark 
 is related to the research for harmonic sums and Mellin transforms
 studied by the physicists Bl\"{u}mlein and Vermaseren \cite{B,V}.
\end{remark}
%
%
\section{The number of the relations}
\label{The-number-of-the-relations}
For $m \ge 0$ we define $V_m$ to be the $\Q$-vector space whose basis
is $I_m$, where $I_0 = \{\phi\}$. In this section, we give the dimension
of the space $u\sigma(V*V) \cap V_m$, which is the number of independent 
relations among MZV's of weight $m$ obtained from Corollary \ref{cor5-20}.
Since $u\sigma$ is an isomorphism of the graded vector space 
$\tV = \bigoplus_{m \ge 0} V_m$, we have
\begin{equation}
 \dim_{\Q}(u\sigma(V*V) \cap V_m) = \dim_{\Q}((V*V) \cap
 V_m). \label{eq5.5-12} 
\end{equation}
Since the bilinear mapping
$* \colon \tilde{V} \times \tilde{V} \to \tilde{V}$
is associative and commutative, $(\tilde{V},*)$ is a commutative $\Q$-algebra
with unit $\phi$. The key to proving Corollary \ref{cor5.5-130} is the fact that
$(\tilde{V},*)$ is the polynomial algebra on the set of
Lyndon words over $\Z_{\ge 1}$, which is due to Malvenuto and Reutenauer~\cite{Ma}.
\par
We start with the definition of the Lyndon words. 
Let $A$ be a set. A
word over $A$ is a finite sequence of elements of $A$. 
We denote the length of a word $u$ by $l(u)$, that is, we put $l(u)=p$
for a word $u=(a_1,\ldots,a_p)$ $(a_i \in A)$.
For two words $u = (a_1, \cdots, a_p)$ and $v = (b_1, \cdots, b_q)$,
we define a word $uv$ by
\begin{displaymath}
 uv = (a_1, \cdots, a_p, b_1, \cdots, b_q).
\end{displaymath}
A word 
$(a_1,a_2,\ldots,a_p)$ is usually written as $a_1a_2\cdots a_p$. 
A lexicographic order on the set of words is given by a total
order $\preceq$ on $A$: For any words $u=a_1 \cdots a_p$, $v=b_1 \cdots b_q$, one
sets $u \preceq v$ if either 
\begin{displaymath}
 p \le q \quad \text{and} \quad a_1=b_1,\,\ldots,\,a_p=b_p 
\end{displaymath}
or
there exists some $1 \le i \le \min\{p,q\}$ such that
\begin{displaymath}
 a_1=b_1,\, \ldots,\, a_{i-1}=b_{i-1},\, a_i \prec b_i.
\end{displaymath}
\begin{definition}
 Let $(A,\preceq)$ be a totally ordered set. A word over $A$ is called a Lyndon word if it is
 strictly smaller than any of its proper right factors, that is, 
$u=a_1 \cdots a_p$ $(a_i \in A)$ is a Lyndon word if 
$u \prec a_i \cdots a_p$ for any $2 \le i \le p$.
\end{definition}
We denote  by $\psi_k(n)$ the number of Lyndon words of length $n$ over
a finite set of cardinality $k$. Then we have 
\begin{equation}
 k^n = \sum_{d|n} d \psi_k(d) \label{eq5.5-20}
\end{equation}
and
\begin{equation}
 n \psi_k(n) = \sum_{d|n} \mu\left(\frac{n}{d}\right) k^d, \label{eq5.5-30}
\end{equation}
where the sums run over the positive divisors of $n$ and 
$\mu$ is the M\"{o}bius function (see \cite[p.~65]{L}).
The equality (\ref{eq5.5-30}) is
obtained from (\ref{eq5.5-20}) and the M\"{o}bius inversion formula. 
\par
In Section \ref{multi-indices}, we defined a multi-index to be a finite
sequence of positive integers. Namely, a multi-index is nothing else but
a word over $\Z_{\ge 1}$. 
Let $\Z_{\ge 1}$ be totally ordered by the usual relation.
We denote the set of Lyndon words over $\Z_{\ge 1}$ by $L$ and the set
of elements in $L$ of weight $m$ by $L_m$.
We define $\Psi(m)$ to be the cardinality of the set $L_m$. 
\begin{proposition}
 \label{prop5.5-70} 
 We have $\Psi(m) = \psi_2(m)$ for any $m \ge 2$.
\end{proposition}
\begin{proof}
 Let $W_m$ be the set of words over $\{x,y\}$ of length $m$ 
 which start with $x$. 
 The elements of $I_m$ are
 in one-to-one correspondence with those of $W_{m}$ by the mapping
 \begin{displaymath}
  (\mu_1,\mu_2,\ldots,\mu_p) \mapsto x \underbrace{y \cdots y}_{\mu_1-1}
  x \underbrace{y \cdots y}_{\mu_2-1} \cdots x \underbrace{y \cdots y}_{\mu_p-1}.
 \end{displaymath}
 If we define a total order on the
 set $\{x,y\}$ by setting $x \prec y$,
 we can easily show that the Lyndon words in $I_m$ are in
 one-to-one correspondence with the Lyndon words in $W_{m}$ under the above
 correspondence. A Lyndon word over $\{x,y\}$ whose length is greater
 than or equal to $2$ necessarily starts with $x$. Hence the
 assertion follows.
\end{proof}
\begin{corollary}
 \label{cor5.5-80} 
 We have
 \begin{displaymath}
  \Psi(m) = \frac{1}{m}\sum_{d|m} \mu\left(\frac{m}{d}\right)2^d
 \end{displaymath}
 for any $m \ge 2$.
\end{corollary}
\begin{proof}
 It follows from (\ref{eq5.5-30}) and Proposition \ref{prop5.5-70}.
\end{proof}
For proofs of the following theorem, see \cite[Section 2]{Hof2} or
\cite[Section 2]{Ma}.
\begin{theorem}
 \label{th5.5-100} 
 $\tV$ is the polynomial algebra over $\Q$ on the set $L$.
\end{theorem}
\begin{corollary}
 \label{cor5.5-130} 
 We have
 \begin{displaymath}
  \dim_{\Q}(u\sigma(V*V) \cap V_m) = 2^{m-1} - \frac{1}{m}\sum_{d|m}\mu\left(\frac{m}{d}\right)2^d
 \end{displaymath}
 for any $m \ge 2$.
\end{corollary}
\begin{proof}
 By Theorem \ref{th5.5-100}, we have
 \begin{displaymath}
  V_m = \Q L_m \oplus \bigl((V*V) \cap V_m \bigr)
 \end{displaymath}
 for any $m \ge 1$. Since $\dim_{\Q} V_m = 2^{m-1}$,
 the assertion follows from Corollary \ref{cor5.5-80} and (\ref{eq5.5-12}).
\end{proof}
%
%
\section{The duality} \label{The-duality} 
In this section, we prove that the assertion of Corollary \ref{cor5-20}
contains the duality. For non-negative integers $p$, $q$ with 
$(p,q) \neq (0,0)$, we define $\calD_{p,q}$ to be the set of 
$2 \times l$ $(l \ge 1)$ matrices $(k_{ij})$ with entries in the
non-negative integers which satisfy $k_{1j} \le 1$ $(1 \le j \le l)$,
$k_{1j} + k_{2j} > 0$ $(1 \le j \le l)$, $\sum_{j=1}^l k_{1j} = p$ and
$\sum_{j=1}^l k_{2j} = q$. For example, $\calD_{1,2}$ consists of 
\begin{multline}
 \begin{pmatrix}
  1 & 0 & 0 \\
  0 & 1 & 1
 \end{pmatrix},\,
 \begin{pmatrix}
  0 & 1 & 0 \\
  1 & 0 & 1
 \end{pmatrix},\,
 \begin{pmatrix}
  0 & 0 & 1 \\
  1 & 1 & 0
 \end{pmatrix}, \\
 \begin{pmatrix}
  1 & 0 \\
  1 & 1
 \end{pmatrix},\,
 \begin{pmatrix}
  0 & 1 \\
  1 & 1
 \end{pmatrix},\,
 \begin{pmatrix}
  1 & 0 \\
  0 & 2
 \end{pmatrix},\,
 \begin{pmatrix}
  0 & 1 \\
  2 & 0
 \end{pmatrix},\,
 \begin{pmatrix}
  1 \\
  2
 \end{pmatrix}. 
 \label{eq6-5}
\end{multline}
We put
\begin{displaymath}
 \calD_{p,q}^0 = \{(k_{ij}) \in \calD_{p,q} \,|\, k_{11} = 0\} \quad
 \text{and} \quad
 \calD_{p,q}^1 = \{(k_{ij}) \in \calD_{p,q} \,|\, k_{11} = 1\}.
\end{displaymath}
For $\bmu = (\muvec)$, $\bnu = (\nuvec) \in \tI$ with 
$(\bmu,\bnu) \neq (\phi,\phi)$, $l(\bmu) = p$ and $l(\bnu) = q$, we
define a mapping
\begin{displaymath}
 \Phi_{p,q;\bmu,\bnu} \colon \calD_{p,q} \to \bigcup_{\bbeta \le \bnu}
 \calH_{\bmu,\bbeta}, \quad 
 \begin{pmatrix}
  k_{11} & \cdots & k_{1l} \\
  k_{21} & \cdots & k_{2l}
 \end{pmatrix}
 \mapsto
 \begin{pmatrix}
  m_{11} & \cdots & m_{1l} \\
  m_{21} & \cdots & m_{2l}
 \end{pmatrix}
\end{displaymath}
by
\begin{multline*}
 (m_{11},m_{12},\ldots,m_{1l}) \\ = (\underbrace{\mu_1 + \mu_2 + \cdots +
 \mu_x}_{k_{11}},\, \underbrace{\mu_{x+1} + \mu_{x+2} + \cdots +
 \mu_y}_{k_{12}},\, \ldots,\, \underbrace{\mu_{z+1} + \mu_{z+2} + \cdots +
 \mu_p}_{k_{1l}})
\end{multline*}
and
\begin{multline*}
 (m_{21},m_{22},\ldots,m_{2l}) \\ = (\underbrace{\nu_1 + \nu_2 + \cdots +
 \nu_{x'}}_{k_{21}},\, \underbrace{\nu_{x'+1} + \nu_{x'+2} + \cdots +
 \nu_{y'}}_{k_{22}},\, \ldots,\, \underbrace{\nu_{z'+1} + \nu_{z'+2} + \cdots +
 \nu_q}_{k_{2l}}),
\end{multline*}
where $\calH_{\bmu,\bbeta}$ is defined in Section \ref{multi-indices}.
For example, $\Phi_{1,2;\bmu,\bnu}$ sends elements of (\ref{eq6-5}) to
\begin{multline}
 \begin{pmatrix}
  \mu_{1} & 0 & 0 \\
  0 & \nu_1 & \nu_2
 \end{pmatrix},\,
 \begin{pmatrix}
  0 & \mu_1 & 0 \\
  \nu_1 & 0 & \nu_2
 \end{pmatrix},\,
 \begin{pmatrix}
  0 & 0 & \mu_1 \\
  \nu_1 & \nu_2 & 0
 \end{pmatrix}, \\
 \begin{pmatrix}
  \mu_1 & 0 \\
  \nu_1 & \nu_2
 \end{pmatrix},\,
 \begin{pmatrix}
  0 & \mu_1 \\
  \nu_1 & \nu_2
 \end{pmatrix},\,
 \begin{pmatrix}
  \mu_1 & 0 \\
  0 & \nu_1 + \nu_2
 \end{pmatrix},\,
 \begin{pmatrix}
  0 & \mu_1 \\
  \nu_1 + \nu_2 & 0
 \end{pmatrix},\,
 \begin{pmatrix}
  \mu_1 \\
  \nu_1 + \nu_2
 \end{pmatrix}, \notag
\end{multline}
respectively. It is easily seen that $\Phi_{p,q;\bmu,\bnu}$ is a bijection.
\begin{proposition} 
 \label{prop6-10} 
 For any $\bmu = (\muvec) \in I$, we have
 \begin{displaymath}
  \sum_{h=0}^p (-1)^h \, \tau(\mu_1,\ldots,\mu_h) * d(\mu_{h+1},\ldots,\mu_p)
  = 0,
 \end{displaymath}
 where $(\mu_1,\ldots,\mu_h) = \phi$ if $h=0$ and
 $(\mu_{h+1},\ldots,\mu_p) = \phi$ if $h=p$.
\end{proposition}
\begin{proof}
 Let $0 \le h \le p$. We denote the composition
 \begin{gather*}
  \mathcal{D}_{h,p-h} \xrightarrow{\Phi_{h,p-h;\tau(\mu_1,\ldots,\mu_h),(\mu_{h+1},\ldots,\mu_p)}} \bigcup_{\bbeta \le
  (\mu_{h+1},\ldots,\mu_p)} \mathcal{H}_{\tau(\mu_1,\ldots,\mu_h),\bbeta}
  \to I\\
  \begin{pmatrix}
   k_{11} & \cdots & k_{1l}\\ k_{21} & \cdots & k_{2l}
  \end{pmatrix}
  \mapsto
  \begin{pmatrix}
   m_{11} & \cdots & m_{1l}\\ m_{21} & \cdots & m_{2l}
  \end{pmatrix}
  \mapsto (m_{11}+m_{21},\ldots,m_{1l}+m_{2l})
 \end{gather*}
 by $\psi_h$. By the definition of the harmonic product $*$ and the
 mapping $d$, we have
 \begin{displaymath}
  \tau(\mu_1,\ldots,\mu_h) * d(\mu_{h+1},\ldots,\mu_p) 
  = \sum_{D \in \mathcal{D}_{h,p-h}} \psi_h(D).
 \end{displaymath}
 For any $0 \le h < p$, the mapping
 \begin{multline*}
  \calD_{h,p-h}^0 \to \calD_{h+1,p-h-1}^1, \\
  D = 
  \begin{pmatrix}
   0 & k_{12} & \cdots & k_{1l} \\
   k_{21} & k_{22} & \cdots & k_{2l}
  \end{pmatrix}
  \mapsto
  D' = 
  \begin{pmatrix}
   1 & k_{12} & \cdots & k_{1l} \\
   k_{21}-1 & k_{22} & \cdots & k_{2l}
  \end{pmatrix}  
 \end{multline*}
 is a bijection and clearly $\psi_h(D) = \psi_{h+1}(D')$. Consequently, we have
 \begin{align*}
  &\sum_{h=0}^p (-1)^h \, \tau(\mu_1,\ldots,\mu_h) *
  d(\mu_{h+1},\ldots,\mu_p)\\
  &= \sum_{h=0}^p (-1)^h \left(\sum_{D \in \calD_{h,p-h}^0} + \sum_{D \in
  \calD_{h,p-h}^1} \right) \, \psi_h(D) \\
  &= \sum_{h=0}^{p-1} (-1)^h \sum_{D \in \calD_{h,p-h}^0} \psi_h(D)
  \quad + \quad \sum_{h=0}^{p-1} (-1)^{h+1} \sum_{D \in \calD_{h+1,p-h-1}^1} \psi_{h+1}(D)\\
  &= 0.
 \end{align*}
 Thus we get the proposition.
\end{proof}
The duality reads
\begin{displaymath}
 (\tau - *)(V) \subset \ker \zeta^{+}.
\end{displaymath}
Therefore Corollary \ref{cor6-10} says that the assertion of Corollary
\ref{cor5-20} contains the duality.
\begin{corollary} 
\label{cor6-10} 
 We have $(\tau - *)(V) \subset u\sigma (V * V)$.
\end{corollary}
\begin{proof}
 Let $\bmu = (\muvec) \in I$. By Proposition \ref{prop6-10}, we have
 \begin{equation}
  -\sum_{h=1}^{p-1} (-1)^h \, \tau (\mu_1,\ldots,\mu_h) *
  d(\mu_{h+1},\ldots,\mu_p) = (\sigma \tau + d)(\bmu). \label{eq6-10} 
 \end{equation}
 By Proposition \ref{eq2-15} and \ref{prop2-20} (v), we have
 $u\sigma(\sigma \tau + d) = (\tau - *)u$.
 Applying $u\sigma$ to both sides of (\ref{eq6-10}), we obtain
 $(\tau - *)u(\bmu) \in u\sigma (V * V)$.
 Thus we have completed the proof.
\end{proof}
%
%
\section{Ohno's relation} \label{Ohno's relation}
In this section, we prove that the assertion of Corollary \ref{cor5-20}
contains Ohno's relation. 
In Section \ref{multi-indices}, we defined $\tilde{V}$ to be the 
$\Q$-vector space whose basis is $\tI$ ($:= \{\phi\} \cup I$).
For $v \in \tV$, we define the $\Q$-linear mapping $m_v$ from $\tV$
to $\tV$ by
\begin{displaymath}
 m_v(w) = v*w \quad (w \in \tV).
\end{displaymath} 
For $\bmu = (\muvec) \in I$, we define the mapping 
$O_{\bmu} \colon \tI \to \tV$ by
\begin{displaymath}
 O_{\bmu}(\bnu) = \sum_{1 \le i_1 < \ldots < i_p \le q}
 (\nu_1,\ldots,\nu_{i_1}+\mu_1,\ldots,\nu_{i_p}+\mu_p,\ldots,\nu_q)
\end{displaymath}
$(\bnu = (\nuvec) \in I)$ and $O_{\bmu}(\phi) = 0$.
We extend the mapping $O_{\bmu}$ linearly on $\tV$.
By setting $O_{\phi} = \id_{\tV}$, we obtain the mapping 
$\tI \to \End_{\Q}(\tV)$, $\bmu \mapsto O_{\bmu}$.
We extend this mapping linearly on $\tV$.
We define $\barO_v = * O_v *$ for any $v \in \tV$.
\begin{proposition} 
\label{prop7-2} 
 For any $\bmu \in \tI$, we have
 $O_{\bmu} \sigma = \sigma O_{\bmu}$ and
 $\barO_{\bmu} \sigma = (-1)^{|\bmu|} \sigma \barO_{\bmu}$.
\end{proposition}
\begin{proof}
 The first assertion is easily seen.
 Let $\bmu = (\mu_1,\ldots,\mu_p) \in I$.
 Let $\bnu \in I$ and  $\bnu^{*} = (\nu^{*}_1,\ldots,\nu^{*}_r)$.
 Then by (\ref{eq2-5}) we have
 \begin{displaymath}
  \sigma \barO_{\bmu}(\bnu)
  = \sum_{1 \le i_1 < \cdots < i_p \le r}(-1)^{|\bmu|+|\bnu|+1-r}*
  (\nu^{*}_1,\ldots,\nu^{*}_{i_1} +
  \mu_1,\ldots,\nu^{*}_{i_p}+\mu_p,\ldots,\nu^{*}_r).
 \end{displaymath}
 Since $|\bnu|+1-r=l(\bnu)$, the second assertion follows.
\end{proof}
\begin{proposition} 
\label{prop7-4} 
 For any $v \in \tV$, we have
 $O_{v} \tau = \tau O_{\tau(v)}$. 
\end{proposition}
\begin{proof}
 This is easily seen.
\end{proof}
\begin{proposition} 
\label{prop7-6} 
 For any $v$, $w \in \tV$, we have $O_v O_w = O_{v*w}$.
\end{proposition}
\begin{proof}
 It is easily seen that $O_{\bmu}O_{\bnu} = O_{\bmu * \bnu}$ 
 for any $\bmu$, $\bnu \in \tI$.
\end{proof}
In the following, if $r=0$, we put
\begin{displaymath}
 (r) = \phi \quad \text{and} \quad (\underbrace{1,\ldots,1}_r) = \phi.
\end{displaymath}
\begin{proposition} 
\label{prop7-10} 
 For any $\bmu$, $\bnu \in I$ and any integer $r \ge 0$, we have
\begin{displaymath}
 O_{u(r)}(\bmu \shp \bnu) = \sum_{\substack{k+l=r\\ k,l
\ge 0}} O_{u(k)}(\bmu) \shp O_{u(l)}(\bnu)
\end{displaymath}
and
\begin{displaymath}
 \barO_{u(r)}(\bmu \dotshp \bnu) = \sum_{\substack{k+l=r\\ k,l
\ge 0}} \barO_{u(k)}(\bmu) \dotshp \barO_{u(l)}(\bnu).
\end{displaymath}
\end{proposition}
\begin{proof}
 The first is easily seen. The second follows from the first and Proposition \ref{prop2-45}.
\end{proof}
\begin{proposition} 
\label{prop7-20} 
 Let $\bmu \in I$. For any $\bnu \in I$, we have
 \begin{displaymath}
   O_{\bmu}((1) \dotshp \bnu) = (1) \dotshp
   O_{\bmu}(\bnu) \quad \text{and} \quad
   \barO_{\bmu}((1) \shp \bnu) = (1) \shp
   \barO_{\bmu}(\bnu).
 \end{displaymath}
\end{proposition}
\begin{proof}
 The first is easily seen. The second
 follows from the first and Proposition \ref{prop2-45}.
\end{proof}
\begin{proposition} 
\label{prop7-30} 
 For any $\bmu \in I$ and any integer $r \ge 0$, we have
 \begin{displaymath}
  O_{\underbrace{1,\ldots,1}_r}(\bmu) 
  = \sum_{\substack{\bnu^1 \shp
  \cdots \shp \bnu^{r+1} = \bmu\\ \bnu^1,\ldots,\bnu^r \in I,\,
  \bnu^{r+1} \in \tI}} \bnu^1 \dotshp \bigl((1) \shp \bnu^2\bigr) \dotshp
  \cdots \dotshp \bigl((1) \shp \bnu^{r+1}\bigr)
 \end{displaymath}
 and
 \begin{displaymath}
  O_{u(r)}(\bmu) 
  = \sum_{\substack{\bnu^1 \shp
  \cdots \shp \bnu^{r+1} = \bmu\\ \bnu^1,\ldots,\bnu^r \in \tI,\,
  \bnu^{r+1} \in I}} \bigl(\bnu^1 \shp (1)\bigr) \dotshp
  \cdots \dotshp \bigl(\bnu^r \shp (1)\bigr) \dotshp
  \bnu^{r+1}.
 \end{displaymath}
\end{proposition}
\begin{proof}
 These are easily seen.
\end{proof}
By Propositions \ref{prop7-30} and \ref{prop2-45}, we obtain
\begin{equation}
 \barO_{u(r)}(\bmu) = \sum_{\substack{\bnu^1 \dotshp
  \cdots \dotshp \bnu^{r+1} = \bmu\\ \bnu^1,\ldots,\bnu^r \in \tI,\,
  \bnu^{r+1} \in I}} (\bnu^1)^{+} \shp 
  \cdots \shp (\bnu^r)^{+} \shp \bnu^{r+1}. \label{eq7-10} 
\end{equation}
Ohno's relation reads
\begin{displaymath}
 O_{u(r)}(t-*)(V) \subset \ker \zeta^{+} \quad (r \ge 0).
\end{displaymath}
The key to proving Corollary \ref{cor7-10},
which says that the assertion of Corollary \ref{cor5-20} contains Ohno's
relation, is the following proposition.
\begin{proposition} 
\label{prop7-40} 
 For any integer $r \ge 0$, we have
 \begin{displaymath}
  u^{-1} m_{\underbrace{1,\ldots,1}_r} u = \sum_{k=0}^r \barO_{u(r-k)}O_{\underbrace{1,\ldots,1}_k}.
 \end{displaymath}
\end{proposition}
\begin{proof}
 It is immediate from Lemmas \ref{lem7-10} and \ref{lem7-20}, which are
 proved later.
\end{proof}
We associate a multi-index $\bmu = (\muvec)$ with a diagram
\begin{displaymath}
 {\setlength{\arraycolsep}{1pt}
 \begin{array}{ccccccccccccccccccc}
           &          &          &          & \downarrow &              &
	   &          &          & \downarrow &            & \downarrow &
	   &          &          &            & \} \text{ arrows} & \\
  \bigcirc & \bigcirc & \cdots   & \bigcirc   &            &\bigcirc    &
  \bigcirc & \cdots   & \bigcirc &            & \cdots     &            &
  \bigcirc & \bigcirc & \cdots     & \bigcirc &            &\:, \\ \cline{1-4}
	   \cline{6-9} \cline{13-16}
           &          & \mu_1      &            &            &            &
	   &
  \mu_2    &          &          &            &            &            &
           & \mu_p    &          &            &            &
 \end{array}}
\end{displaymath}
where the number of arrows is $p-1$. For example, we associate
multi-indices $(2,2)$ and $(1,2,1)$ with diagrams
\begin{displaymath}
{\setlength{\arraycolsep}{1pt}
 \begin{array}{ccccccccccccccc}
                         & \phantom{\downarrow}   &
			 & \downarrow		  &
			 & \phantom{\downarrow}   &          
                         &	                  &
			 & \downarrow             & 
                         & \phantom{\downarrow}   &
                         & \downarrow             & \\  
  \bigcirc               &                        & \bigcirc
                         &                        & \bigcirc
                         &                        & 
  \bigcirc               & \quad \text{and} \quad & \bigcirc
                         &                        & \bigcirc
			 &                        & \bigcirc
                         &                        & \bigcirc \,\,.\\
 \end{array}}
\end{displaymath}
Then multi-indices are in one-to-one
correspondence with these diagrams.
Let $n_1$, \ldots, $n_r$ ($r \geq 1$) be non-negative integers with
$n_1 + \cdots + n_r = |\bmu|$. We define $\bnu^1$, \ldots, 
$\bnu^r \in \tI$ by\\[0.5mm]
\begin{displaymath}
 \hspace{-7em} \overbrace{\hspace{4em}}^{\bnu^1} \hspace{0.5em} \overbrace{\hspace{4.5em}}^{\bnu^2} 
 \hspace{3em} \overbrace{\hspace{4em}}^{\bnu^r} 
\end{displaymath}
\vspace{-1.5em}
\begin{displaymath} 
 {\setlength{\arraycolsep}{1.5pt}
 \begin{array}{ccccccccccccccc}
           &          &            &             &          &
           &          &            &             &          &
           &          &            & \} \text{ arrows of } \bmu & \\ 
  \bigcirc & \cdots   & \bigcirc   & \phantom{a} & \bigcirc & 
  \cdots   & \bigcirc & \phantom{a}& \cdots      & \phantom{a} & 
  \bigcirc & \cdots   & \bigcirc   &             & \:, \\ \cline{1-3} \cline{5-7}
  \cline{11-13}
           & n_1      &            &             &          & 
    n_2    &          &            &             &          &
           & n_r      &            &             &
 \end{array}}
\end{displaymath}\\[0.5mm]
where the weight of $\bnu^i$ is $n_i$. We call
$(\bnu^1,\ldots,\bnu^r)$ the $(n_1,\ldots,n_r)$-partition of $\bmu$
and write $(\bnu^1,\ldots,\bnu^r) \,||\, \bmu$.
For example, the $(2,1,3)$-partition and the $(3,0,3)$-partition of 
$(2,2,2)$ are $((2),\,(1),\,(1,2))$ and $((2,1),\,\phi,\,(1,2))$ by the
following diagrams
\begin{displaymath}
 {\setlength{\arraycolsep}{1pt}
 \begin{array}{cccccccccccccccccccccccc}
           & \phantom{\downarrow} & & \downarrow & &
  \phantom{\downarrow} & & \downarrow & & \phantom{\downarrow} & & & &
  \phantom{\downarrow} & & \downarrow & & \phantom{\downarrow} & &
  \downarrow & & \phantom{\downarrow} & & \\
  \bigcirc & & \bigcirc & & \bigcirc & & \bigcirc & &\bigcirc & &
  \bigcirc & \quad \text{and} \quad & \bigcirc & & \bigcirc & & \bigcirc 
           & & \bigcirc & & \bigcirc & & \bigcirc & \,\,. \\
  \cline{1-3} \cline{5-5} \cline{7-11} \cline{13-17} \cline{19-23}
  & 2 & & & 1 & & & & 3 & & & & & & 3 & & & & & & 3 & & &
 \end{array}}
\end{displaymath}
\begin{lemma} 
\label{lem7-10} 
 For any $\bmu \in I$ and any integer $r \ge 0$, we have
 \begin{displaymath}
  u^{-1} m_{\underbrace{1,\ldots,1}_r} u(\bmu) =
  \sum_{\substack{(\bnu^1,\ldots,\bnu^{r+1}) \,||\, \bmu\\
  \bnu^i \in \tI}} (\bnu^1)^{+} \shp \cdots \shp (\bnu^r)^{+}
  \shp \bnu^{r+1}.
 \end{displaymath}
\end{lemma}
\begin{proof}
 We define a mapping $f \colon \tI \to V$ by $f(\phi) = (1)$ and
 $f(\balpha) = \balpha \shp (1) + \balpha \dotshp (1)$ $(\balpha \in I)$.
 Then we have
 \begin{displaymath}
  m_{\underbrace{1,\ldots,1}_r}(\bmu) = \sum_{\substack{\balpha^1 \shp
  \cdots \shp \balpha^{r+1} = \bmu\\ \balpha^i \in \tI}}
  f(\balpha^1) \shp \cdots \shp f(\balpha^r) \shp \balpha^{r+1}
 \end{displaymath}
 and therefore
 \begin{equation}
  m_{\underbrace{1,\ldots,1}_r}u(\bmu) 
  = \sum_{\substack{\balpha^1 \shp \cdots \shp \balpha^{r+1} \ge
  \bmu\\ \balpha^i \in \tI}}
  f(\balpha^1) \shp \cdots \shp f(\balpha^r) \shp
  \balpha^{r+1}. \label{eq7-20} 
 \end{equation}
 We construct a mapping
 \begin{multline*}
  \bigl\{(\balpha^1,\ldots,\balpha^{r+1}) \in \tI^{r+1} \,|\, \balpha^1 \shp \cdots
  \shp \balpha^{r+1} \ge \bmu \bigr\} \to\\
  \biggl\{\bigl((\balpha^1,\ldots,\balpha^{r+1}),\,(\bnu^1,\ldots,\bnu^{r+1})\bigr)
  \in \tI^{r+1} \times \tI^{r+1} \biggm|
  \begin{array}{c}
   (\bnu^1,\cdots,\bnu^{r+1}) \,||\, \bmu,\\ 
   \balpha^i \ge \bnu^i \, (1 \le i \le r+1) 
  \end{array}
  \biggr\}, \\
  (\balpha^1,\ldots,\balpha^{r+1}) \mapsto \bigl((\balpha^1,\ldots,\balpha^{r+1}),\,(\bnu^1,\ldots,\bnu^{r+1})\bigr)
 \end{multline*}
 by defining
 $(\bnu^1,\ldots,\bnu^{r+1})$ to be the
 $(|\balpha^1|,\ldots,|\balpha^{r+1}|)$-partition of $\bmu$.
 The mapping
 \begin{displaymath}
  \bigl((\balpha^1,\ldots,\balpha^{r+1}),\,(\bnu^1,\ldots,\bnu^{r+1})\bigr)
  \mapsto (\balpha^1,\ldots,\balpha^{r+1})
 \end{displaymath} 
 is its inverse and thereby
 this mapping is a bijection. Hence, the right-hand side of
 (\ref{eq7-20}) is equal to
 \begin{align*}
  &\sum_{\substack{(\bnu^1,\ldots,\bnu^{r+1}) \,||\, \bmu\\
  \bnu^i \in \tI}}
  \sum_{\substack{\balpha^1 \ge \bnu^1\\ \cdots\\ \balpha^{r+1} \ge
  \bnu^{r+1}}}
  f(\balpha^1) \shp \cdots \shp f(\balpha^r) \shp \balpha^{r+1}\\
  &= \sum_{\substack{(\bnu^1,\ldots,\bnu^{r+1}) \,||\, \bmu\\
  \bnu^i \in \tI}}
  \left\{\sum_{\balpha^1 \ge \bnu^1} f(\balpha^1)\right\} \shp
  \cdots \shp 
  \left\{\sum_{\balpha^r \ge \bnu^r} f(\balpha^r)\right\} \shp \left\{\sum_{\balpha^{r+1} \ge \bnu^{r+1}}
  \balpha^{r+1} \right\} \\
  &= \sum_{\substack{(\bnu^1,\ldots,\bnu^{r+1}) \,||\, \bmu\\
  \bnu^i \in \tI}}
  u\bigl((\bnu^1)^{+}\bigr) \shp \cdots \shp u\bigl((\bnu^r)^{+}\bigr) \shp u(\bnu^{r+1}).
 \end{align*}
 By Proposition \ref{prop2-50}, we obtain the assertion.
\end{proof}
\begin{remark} 
\label{rem7-10}
 We have $u^{-1} = \sigma u \sigma$ by Proposition \ref{prop2-20} (iii).
 According to Corollary \ref{cor5-20}, we obtain
 \begin{displaymath}
  \sigma \sum_{\substack{(\bnu^1,\ldots,\bnu^{r+1}) \,||\,
  \bmu\\ \bnu^i \in \tI}}
  (\bnu^1)^{+} \shp \cdots \shp (\bnu^r)^{+} \shp \bnu^{r+1}
  \in \ker \zeta^{+}
 \end{displaymath}
 for any $\bmu \in I$ and any integer $r \ge 1$.
 For example, for $\bmu = (2,1)$ and $r=2$,
 all of partitions $(\bnu^1,\bnu^2)$ of $\bmu$ are
 $\bigl((2),(1)\bigr)$, $\bigl((1),(1,1)\bigr)$, 
 $\bigl((2,1),\phi \bigr)$, $\bigl(\phi,(2,1)\bigr)$
 by the diagrams
 \begin{displaymath}
  {\setlength{\arraycolsep}{1pt}
  \begin{array}{cccccccccccccccccc}
   & \phantom{\downarrow} & & \downarrow & & \quad & &
    \phantom{\downarrow} & & \downarrow & & & &  \phantom{\downarrow} &
    & \downarrow & & \\
   \bigcirc & & \bigcirc & & \bigcirc & \, ,\quad & \bigcirc & & \bigcirc &
    & \bigcirc & \quad \text{and} \quad & \bigcirc & & \bigcirc & &
    \bigcirc &. \\ \cline{1-3} \cline{5-5} \cline{7-7} \cline{9-11} \cline{13-17}
  \end{array}}
 \end{displaymath}
 Since $(2)^{+} \shp (1) = (3,1)$, $(1)^{+} \shp (1,1) = (2,1,1)$,
 $(2,1)^{+} \shp \phi = (2,2)$ and $\phi^{+} \shp (2,1) = (1,2,1)$,
 we obtain $(3,1) - (2,1,1) + (2,2) - (1,2,1) \in \ker \zeta^{+}$.
\end{remark}
\begin{lemma} 
\label{lem7-20} 
 For any $\bmu \in I$ and any integer $r \ge 0$, we have
 \begin{displaymath}
  \sum_{k=0}^r \barO_{u(r-k)} O_{\underbrace{1,\ldots,1}_{k}}(\bmu) =
  \sum_{\substack{(\bnu^1,\ldots,\bnu^{r+1}) \,||\, \bmu\\ \bnu^i \in
  \tI}} (\bnu^1)^{+} \shp \cdots \shp (\bnu^r)^{+} \shp \bnu^{r+1}.
 \end{displaymath}
\end{lemma}
\begin{proof}
 By Proposition \ref{prop7-30}, we have
 \begin{displaymath}
  O_{\underbrace{1,\ldots,1}_{k-1}}(\bmu) =
  \sum_{\substack{\balpha^1 \shp \cdots \shp \balpha^k = \bmu\\
  \balpha^1,\ldots,\balpha^{k-1} \in I,\,\balpha^k \in \tI}}
  \balpha^1 \dotshp \bigl((1) \shp \balpha^2\bigr) \dotshp \cdots \dotshp
  \bigl((1) \shp \balpha^k \bigr).
 \end{displaymath}
 Therefore, by Propositions \ref{prop7-10} and \ref{prop7-20},
 we have
 \begin{multline*}
  \barO_{u(r-k+1)} O_{\underbrace{1,\ldots,1}_{k-1}}(\bmu)
  =
  \sum_{\substack{\balpha^1 \shp \cdots \shp \balpha^k = \bmu\\
  \balpha^1,\ldots,\balpha^{k-1} \in I,\,\balpha^k \in \tI}}
  \sum_{\substack{m_1 + \cdots + m_k = r+1\\ m_i \ge 1}} \\
  \barO_{u(m_1-1)}(\balpha^1) \dotshp \bigl((1) \shp
  \barO_{u(m_2-1)}(\balpha^2)\bigr) \dotshp \cdots 
  \dotshp
  \bigl((1) \shp \barO_{u(m_k-1)}(\balpha^k)\bigr),
 \end{multline*}
 where we put
 \begin{displaymath}
  \barO_{u(m_k-1)}(\balpha^k) = (\underbrace{1,\ldots,1}_{m_k-1})
 \end{displaymath}
 if $\balpha^k = \phi$. By (\ref{eq7-10}), this is equal to
 \begin{align*}
  &\sum_{\substack{\balpha^1 \shp \cdots \shp \balpha^k = \bmu\\
  \balpha^1,\ldots,\balpha^{k-1} \in I,\,\balpha^k \in \tI}}
  \sum_{\substack{m_1 + \cdots + m_k = r+1\\ m_i \ge 1}} \\
  &\hspace{3em}\left\{\sum_{\substack{\blambda^{11} \dotshp \cdots
  \dotshp \blambda^{1m_1} = \balpha^1\\
  \blambda^{11},\ldots,\blambda^{1m_1-1} \in \tI,\, \blambda^{1m_1}
  \in I}}
  (\blambda^{11})^{+} \shp \cdots \shp (\blambda^{1m_1-1})^{+} \shp
  \blambda^{1m_1} \right\} \\
  &\hspace{2em} \dotshp \left\{(1) \shp \sum_{\substack{\blambda^{21} \dotshp \cdots
  \dotshp \blambda^{2m_2} = \balpha^2\\
  \blambda^{21},\ldots,\blambda^{2m_2-1} \in \tI,\, \blambda^{2m_2}
  \in I}}
  (\blambda^{21})^{+} \shp \cdots \shp (\blambda^{2m_2-1})^{+} \shp
  \blambda^{2m_2} \right\} \\
  &\hspace{2em}\dotshp \cdots \\
  &\hspace{2em} \dotshp \left\{(1) \shp \sum_{\substack{\blambda^{k1} \dotshp \cdots
  \dotshp \blambda^{km_k} = \balpha^k\\
  \blambda^{k1},\ldots,\blambda^{km_k} \in \tI\\ \balpha^k \neq \phi
  \,\Rightarrow \, \blambda^{km_k} \neq \phi}}
  (\blambda^{k1})^{+} \shp \cdots \shp (\blambda^{km_k-1})^{+} \shp
  \blambda^{km_k} \right\} \\
  &=\sum_{\substack{\balpha^1 \shp \cdots \shp \balpha^k = \bmu\\
  \balpha^1,\ldots,\balpha^{k-1} \in I\\ \balpha^k \in \tI}} \,\,
  \sum_{\substack{m_1 + \cdots + m_k = r+1\\ m_i \ge 1}} \,\,
  \sum_{\substack{\blambda^{11} \dotshp \cdots
  \dotshp \blambda^{1m_1} = \balpha^1\\
  \blambda^{11},\ldots,\blambda^{1m_1-1} \in \tI\\ \blambda^{1m_1}
  \in I}} \cdots 
  \sum_{\substack{\blambda^{k1} \dotshp \cdots
  \dotshp \blambda^{km_k} = \balpha^k\\
  \blambda^{k1},\ldots,\blambda^{km_k} \in \tI\\ \balpha^k \neq \phi
  \,\Rightarrow \, \blambda^{km_k} \neq \phi}}\\
  &\hspace{3em} (\blambda^{11})^{+} \shp \cdots \shp 
  (\blambda^{1m_1})^{+}
  \shp \cdots \shp 
  (\blambda^{k1})^{+} \shp \cdots \shp (\blambda^{km_k-1})^{+} \shp
  \blambda^{km_k}.
 \end{align*}
 We define a mapping
 \begin{multline*}
  \Bigl\{(\bnu^1,\ldots,\bnu^{r+1}) \in \tI^{r+1} \bigm|
  (\bnu^1,\ldots,\bnu^{r+1}) \,||\, \bmu \Bigr\} \to \\
  \shoveleft{\Bigl\{\bigl((\blambda^{11},\ldots,\blambda^{1m_1}),\ldots,(\blambda^{k1},\ldots,\blambda^{km_k})\bigr)
  \Bigm| \Bigr.}\\
  \left.
  \begin{array}{c}
  1 \le k \le r+1,\,
  m_1 + \cdots + m_k = r+1,\,m_i \ge 1,\\ 
  (\blambda^{11} \dotshp \cdots \dotshp \blambda^{1m_1}) \shp \cdots
   \shp (\blambda^{k1} \dotshp \cdots \dotshp \blambda^{km_k}) = \bmu, \\
  \blambda^{ij} \in \tI,\, \blambda^{1m_1},\ldots,\blambda^{k-1
  m_{k-1}} \neq \phi,\, 
  \blambda^{km_k} = \phi \,\Rightarrow \,
  \blambda^{k1}=\cdots=\blambda^{km_k}=\phi
  \end{array}
  \right\}, \\
  (\bnu^1,\ldots,\bnu^{r+1}) \mapsto \bigl((\blambda^{11},\ldots,\blambda^{1m_1}),\ldots,(\blambda^{k1},\ldots,\blambda^{km_k})\bigr)
 \end{multline*}
 by
 \begin{displaymath}
  (\blambda^{11},\ldots,\blambda^{1m_1},\ldots,\blambda^{k1},\ldots,\blambda^{km_k})
  = (\bnu^1,\ldots,\bnu^{r+1})
 \end{displaymath}
 and
 \begin{multline*}
  \{m_1,\, m_1 + m_2,\, \ldots,\, m_1 + m_2 + \cdots + m_k \} \\
  = \left\{1 \le i \le r+1 \biggm| 
  \begin{array}{l}
   |\bnu^1| + \cdots + |\bnu^i| = \mu_1 + \cdots
   + \mu_j \\ 
   \text{ for some } \, 1 \le j \le l(\bmu) 
   \text{ and } \, \bnu^i \neq \phi \,\text{ if }\, i \neq r+1
  \end{array}\right\}. 
 \end{multline*}
 The following diagram illustrates the definition of this mapping.
\begin{align*}
 &\overbrace{\hspace{3.5em}}^{\blambda^{11}} \hspace{3em}
 \overbrace{\hspace{4em}}^{\blambda^{1m_1}} \hspace{0.8em}
 \overbrace{\hspace{4.8em}}^{\blambda^{21}} \hspace{3.2em}
 \overbrace{\hspace{4em}}^{\blambda^{2m_2}} \\[-2mm]
 &{\setlength{\arraycolsep}{0pt}
 \begin{array}{ccccccccccccccccccccc}
  & & & \times & & \times & & & & \downarrow & & & & \times & & \times &
   & & & \downarrow & \\
  \bigcirc & \cdots & \bigcirc & &\cdots & & \bigcirc & \cdots &
   \bigcirc & & \bigcirc & \cdots & \bigcirc & & \cdots & & \bigcirc &
   \cdots & \bigcirc & & \cdots \\
  \cline{1-3} \cline{7-9} \cline{11-13} \cline{17-19}
  & |\bnu^1| & & & & & & |\bnu^{M_1}| & & & & |\bnu^{M_1+1}| & & & & & &
  |\bnu^{M_2}| & & & 
 \end{array}}\\[3mm]
 &\hspace{12em} \overbrace{\hspace{5.7em}}^{\blambda^{k1}} \hspace{3.2em} \overbrace{\hspace{4.2em}}^{\blambda{km_k}}\\[-2mm]
 &\hspace{10em}{\setlength{\arraycolsep}{0pt}
 \begin{array}{cccccccccccc}
  &\downarrow & & & & \times & & \times & & & & \} \text{ arrows of } \bmu \\
  \cdots && \bigcirc & \cdots & \bigcirc & & \cdots & & \bigcirc & \cdots &
   \bigcirc & \\
  \cline{3-5} \cline{9-11}
  && & |\bnu^{M_{k-1}+1}| & & & & & & |\bnu^{M_k}| & & \\
 \end{array}}
\end{align*}
\begin{center}
 ($M_i = m_1 + \cdots + m_i$. There exist no arrows of $\bmu$ at the
 positions of $\times$.)
\end{center} 
The inverse is given by
\begin{align*}
 &\bigl((\blambda^{11},\ldots,\blambda^{1m_1}),\ldots,(\blambda^{k1},\ldots,\blambda^{km_k})\bigr) \\
 &\hspace{10em}\mapsto
 (\blambda^{11} \dotshp \cdots \dotshp \blambda^{1m_1}) \shp \cdots \shp (\blambda^{k1}
 \dotshp \cdots \dotshp \blambda^{km_k}),
\end{align*}
and therefore the mapping is a bijection. As a result, we obtain
\begin{displaymath}
 \sum_{k=1}^{r+1}
 \barO_{u(r-k+1)}O_{\underbrace{1,\ldots,1}_{k-1}}(\bmu)
 = 
 \sum_{\substack{(\bnu^1,\ldots,\bnu^{r+1}) \,||\, \bmu\\ \bnu^i \in
 \tI}}
 (\bnu^1)^{+} \shp \cdots \shp (\bnu^r)^{+} \shp \bnu^{r+1}.
\end{displaymath}
We have thus completed the proof.
\end{proof}
\begin{proposition} 
\label{prop7-50} 
 For any integer $r \ge 0$, we have
 \begin{displaymath}
  \sum_{k=0}^r (-1)^k u^{-1}m_{\underbrace{1,\ldots,1}_{r-k}}u \,O_{u(k)}
  = \barO_{u(r)}.
 \end{displaymath}
\end{proposition}
\begin{proof}
 By Proposition \ref{prop7-40}, the left-hand side is equal to
 \begin{equation}
  \sum_{r \ge l \ge k \ge 0} (-1)^k \barO_{u(r-l)}  
  O_{\underbrace{1,\ldots,1}_{l-k}} O_{u(k)}. \label{eq7-30} 
 \end{equation}
 By Propositions \ref{prop7-6} and \ref{prop6-10}, if $l > 0$, we have
 \begin{displaymath}
  \sum_{k=0}^l (-1)^{k} O_{\underbrace{1,\ldots,1}_{l-k}}O_{u(k)}
  = \sum_{k=0}^l (-1)^{k} O_{(\underbrace{1,\ldots,1}_{l-k}) *
  d(\underbrace{1,\ldots,1}_k)} = 0. 
 \end{displaymath}
 Therefore (\ref{eq7-30}) is equal to $\barO_{u(r)}$.
\end{proof}
Since $u^{-1} = \sigma u \sigma$, by Proposition \ref{prop7-50} we have
\begin{displaymath}
 \sigma \bigl(\barO_{u(r)} - (-1)^r O_{u(r)} \bigr) = \sum_{k=0}^{r-1}(-1)^k u\sigma
 m_{\underbrace{1,\ldots,1}_{r-k}}u
  \,O_{u(k)}
\end{displaymath}
and therefore
\begin{equation}
 \sigma \bigl(\barO_{u(r)} - (-1)^r O_{u(r)} \bigr)(V) \subset u\sigma
  (V*V). \label{eq7-50} 
\end{equation}
\begin{corollary} 
\label{cor7-10} 
 For any integer $r \ge 0$, we have
 $O_{u(r)}(\tau - *)(V) \subset u\sigma (V*V)$.
\end{corollary}
\begin{proof}
 Since 
 \begin{displaymath}
  \sigma \bigl((-1)^r \barO_{u(r)} - O_{u(r)}\bigr)(V) =
  \bigl(\barO_{u(r)} - O_{u(r)} \bigr)(V)
 \end{displaymath}
 by Proposition \ref{prop7-2}, we have
 \begin{equation}
  \bigl(\barO_{u(r)} - O_{u(r)} \bigr)(V) \subset u\sigma(V*V)
   \label{eq7-60} 
 \end{equation}
 by (\ref{eq7-50}). Since $O_{u(r)}\tau = \tau O_{u(r)}$ by Proposition \ref{prop7-4}, we have
 \begin{displaymath}
  O_{u(r)}\tau - O_{u(r)}* = (\tau - *)O_{u(r)} + \bigl(\barO_{u(r)} - O_{u(r)} \bigr)*.
 \end{displaymath}
 This equality together with (\ref{eq7-60}) and Corollary \ref{cor6-10} implies
 the assertion of Corollary \ref{cor7-10}.
\end{proof}
\noindent{\textbf{Acknowledgments.}}
The author would like to express his gratitude to Professor Masanobu Kaneko and Professor Yasuo Ohno for valuable comments and suggestions.
Especially Professor Ohno recommended the author to calculate the numerical values of the dimensions of our relations
and Professor Kaneko suggested the dimension formula as is stated in the introduction.
The author would also like to thank Professor Yoshio Tanigawa for
helpful advices.
Finally the author would like to thank the anonymous referee for
a number of valuable comments to improve the presentation of the paper.

\end{document}